\documentclass{amsart}
\usepackage{amsmath,amssymb,amsxtra,enumerate,graphics}
\usepackage[mathscr]{eucal}

\newtheorem{theorem}{Theorem}[section]
\newtheorem{lemma}[theorem]{Lemma}
\newtheorem{proposition}[theorem]{Proposition}
\newtheorem{corollary}[theorem]{Corollary}

\newcommand{\Psibar}{\overline\Psi}

\newcommand{\rss}{\text{rss}}

\newcommand{\tbf}{\textbf}
\newcommand{\addots}{\reflectbox{$\ddots$}}  

\newcommand{\al}{{\alpha}}
\newcommand{\be}{\beta}
\newcommand{\ga}{\gamma}

\newcommand{\del}{\delta}

\newcommand{\ep}{\varepsilon}

\newcommand{\Z}{\mathbb{Z}}

\newcommand{\sF}{\mathscr{F}}  

\newcommand{\kbar}{\bar{k}}

\newcommand{\wdot}{{\dot{w}}}

\newcommand{\ds}{\oplus}

\newcommand{\bigds}{\bigoplus}

\newcommand{\goesto}{\mapsto}

\newcommand{\bmat}{\begin{pmatrix}}
\newcommand{\emat}{\end{pmatrix}}

\newcommand{\GL}{\text{\rm GL}}

\newcommand{\gl}[2]{\text{\rm GL}_{#1}\!\left(#2\right)}
\renewcommand{\sl}[2]{\text{\rm SL}_{#1}\!\left(#2\right)}

\newcommand{\offd}[2]{{#1\!\left(#2\right)}}

\DeclareMathOperator{\ad}{ad}

\DeclareMathOperator{\Aut}{Aut}

\DeclareMathOperator{\Sym}{Sym}

\DeclareMathOperator{\height}{ht}

\DeclareMathOperator{\vol}{{vol}}
\DeclareMathOperator{\CA}{{A}}
\DeclareMathOperator{\CB}{{B}}
\DeclareMathOperator{\CC}{{C}}
\DeclareMathOperator{\CD}{{D}}
\DeclareMathOperator{\CE}{{E}}
\DeclareMathOperator{\CF}{{F}}
\DeclareMathOperator{\CG}{{G}}

\DeclareMathOperator{\SO}{{SO}}

\newcounter{marks}
\newcounter{question}

\begin{document}
\title[Algorithm for Lang's Theorem]{Algorithm for Lang's Theorem}
\author{Arjeh M.\ Cohen}
\author{Scott H. Murray}
\date{\today}
\thanks{The authors would like to thank Sergei Haller, Anthony Henderson,
William M.~Kantor, Gus Lehrer, T.~A. Springer, and D.~E. Taylor
for useful discussions on this topic.  The authors also thank the
Magma project at the University of Sydney, where some of the work was carried
out.}

\begin{abstract}
  We give an efficient algorithm for Lang's Theorem in split connected
  reductive groups defined over finite fields of characteristic greater than 3.
  This algorithm can be used to construct many important structures in
  finite groups of Lie type.
  We use an algorithm for computing a Chevalley basis for a split
  reductive Lie algebra, which is of independent interest.
\end{abstract}

\maketitle

\section{Introduction}\label{S-intro}
A finite group of Lie type can be described as the rational points
of a connected reductive algebraic group over a finite field.
Given a structure in the algebraic group, such as a conjugacy class or a maximal
torus, we want to find the corresponding structures in the finite group
of Lie type.
This can often be achieved with Lang's Theorem.
We provide a computationally efficient algorithm for Lang's Theorem in
split connected reductive groups.
Our algorithm is randomised but guaranteed to return a correct answer,
ie, it is Las Vegas in the sense of \cite{Babai97}.
Glasby and Howlett \cite{GlasbyHowlett97} have already solved this
problem in a special case;  our algorithm is inspired by their work and the
proof of Lang's Theorem given in \cite{Muller03}.

Throughout this paper $k$ is a finite field of size $q$ and characteristic $p$,
and $k_r$ is the unique degree $r$ extension of $k$ in the
algebraic closure $\kbar$.
The affine space of dimension $N$ can be identified with $\kbar^N$.
An \emph{affine variety} $X$ is a subset of $\kbar^N$ that consists of
the zeroes of a collection of polynomials.
The variety is \emph{defined over $k$} if it is closed under the action
of the map $F:\kbar^N\to\kbar^N$ that takes the $q$th power
of each component.
The restriction of $F$ to $X$ is called the
\emph{(standard) Frobenius endomorphism} of $X$.
The set of \emph{rational points} of $X$ over $k_r$, denoted by $\offd{X}{k_r}$,
consists of those elements of $X$ fixed by $F^{r}$.
A \emph{nonstandard Frobenius endomorphism} is a morphism $F':X\to X$ such that
$(F')^s=F^s$ for some positive integer $s$.
The elements of $X$ fixed by $F'$ are the rational points of a 
\emph{$k$-form} of $X$. 
In this paper, Frobenius endomorphisms are standard unless otherwise stated.

A \emph{linear algebraic group} is an affine variety with group multiplication and
inversion given by rational functions.
See, for example, \cite{Springer98Book} for more details including the
definitions of reductive and connected groups.
Every linear algebraic group contains a maximal connected subgroup $G^\circ$,
the component of the identity.
This subgroup is normal and $G/G^\circ$ is finite, so for many purposes it
suffices to study connected groups.
An important result on  linear algebraic groups over finite
fields is:
\begin{theorem}[Lang's Theorem]
If $G$ is a connected linear algebraic group defined over the finite field $k$ with Frobenius
map $F$, then the map
$$
  G\to G,\quad a\goesto a^{-F}a
$$
is onto.
\end{theorem}
\noindent
This is equivalent to the statement that the first Galois cohomology 
of $G$
 is trivial.

In this paper, we give an algorithm for Lang's Theorem in $k$-split connected
reductive groups 
described by the
Steinberg presentation as in \cite{CohenMurrayTaylor04}.
In particular, the root datum, and hence the Cartan type, of $G$ is known.
Reductive groups are likely to be the critical case, since the problem
for an arbitrary connected linear algebraic group  could be solved by working down
a composition series (see Section~\ref{S-evect}) and all simple connected 
groups are reductive.
Our main result is:
\begin{theorem}\label{T-time}
Let $k$ be a finite field of characteristic greater than $3$.
Let $G$ be a $k$-split connected reductive linear algebraic group.
Let $c$ be in $G(k_r)$, and suppose we are given $s$, the order of
$c^{F^{r-1}}\cdots c^Fc$.
Then we can find $a\in G(k_{rs})$ such that $c=a^{-F}a$ in 
Las Vegas time $O(n^9r^2s^2\log^2(n)\log^2(q))$
where $n$ is the reductive rank of $G$.
\end{theorem}
\noindent
We can improve significantly on this result for the classical groups:
\begin{theorem}\label{T-classtime}
Let $G$ be a $k$-split simple connected classical group defined over $k$.  
Let $c$ be in $G(k_r)$ and suppose we are given $s$, the order of
$c^{F^{r-1}}\cdots c^Fc$.
Then we can find $a\in G(k_{rs})$ such that $c=a^{-F}a$ in
Las Vegas time $O(n^5r^2s^2\log^2(q))$ where $n$ is the reductive rank of $G$. 
\end{theorem}

The parameter $s$ measures the size of the field extension required, as explained
in Section~\ref{S-min}. 
In Section~\ref{S-evect}, we use the concept of $F$-eigenvectors to reduce to a
problem involving forms of $G$-modules.
A solution to this problem and a proof of Theorem~\ref{T-classtime}
is given in Section~\ref{S-transf}.
This solution uses the algorithm for computing a standard Chevalley basis in the
Lie algebra of $G$ described in Section~\ref{S-csa}.
The running time of this algorithm is analysed in Section~\ref{S-prob}, leading
to a proof of Theorem~\ref{T-time}.

\section{Minimum field degree}\label{S-min}
Computation in large finite fields is a challenging problem
(see, for example, \cite{LidlNiederreiter97Book}).
So we start with an easy result giving the
size of the field extension needed for Lang's Theorem.
We define the \emph{minimum field degree} of $g\in G$ as the
smallest $r$ such that $g^{F^r}=g$.
Note that $g$ has minimum field degree $r$ if, and only if, $k_r$ is the
smallest extension of $k$ such that $g$ is in $\offd{G}{k_r}$.

\begin{proposition}\label{P-min}
Let $G$ be a connected linear algebraic group defined over $k$.
Let $c$ be an element of $G$ with minimum field degree $r$ and
let $s$ be the order of $c^{F^{r-1}}\cdots c^Fc$.
If $c=a^{-F}a$ for some $a$ in $G$, then
the minimum field degree of $a$ is~$rs$.
\end{proposition}
\begin{proof}
Let $m$ be the minimum field degree of $a$.
Clearly $k_r$ is a subfield of $k_m$, so $r$ is a divisor of $m$, say $ru=m$.
Since $c^{F^r}=c$, we have
$$
  \left(c^{F^{r-1}}\cdots c^Fc\right)^u = c^{F^{m-1}}\cdots c^Fc
    = a^{-F^m}a^{F^{m-1}}\cdots a^{-F^2}a^{F}a^{-F}a= a^{-F^m}a.
$$
Hence $a^{F^m}=a$ if, and only if, $u$ is a multiple of $s$.
\end{proof}

The most important consequence of this proposition is that the minimum field degree
is independent of the particular choice of $a$
and can be computed beforehand.
In all our timings we consider $s$, the order of $c^{F^{r-1}}\cdots c^Fc$,
to be an input of our algorithm.
While it is straightforward to compute $s$, no polynomial time algorithm is
known.
The best known method for computing $s$ is to convert from the Steinberg
presentation of $G$ to a faithful representation \cite{CohenMurrayTaylor04} and then
compute the order of the corresponding matrix using the algorithm of
\cite{CellerLeedham-Green97}.
If the representation has degree $d$,
this takes Las Vegas time $O(d^3\log(q)\log\log(q^d))$ plus the time required to
factor a collection of integers of the form $q^{d_i}-1$ with $\sum_i d_i\le d$.

Suppose now that $G$ is a $k$-split reductive group with reductive rank $n$ and
semisimple rank $\ell$.
The element $c$, which is the input to our algorithm, has size
$O((n+\ell^2)r\log(q))$; while the element $a$, which is the output,
has size $O((n+\ell^2)rs\log(q))$.
Since $s$  need not be bounded by a polynomial in $n$, $\ell$, $r$, and $\log(q)$,
there is no algorithm for Lang's Theorem that is polynomial in the
size of the input.
The best we can hope for is an algorithm which is polynomial
in the size of the output.

\section{Twisted eigenvectors}\label{S-evect}
We can now give an outline of our main algorithm.
Let $G=\offd{G}{\kbar}$ be a connected linear algebraic group defined over $k$.
Suppose that $V$ is a $G$-module of dimension $d$ defined over $k$, so that
$F$ acts on $V=\kbar^d$ by taking the $q$th power of each component.
We say that $v\in V$ is an \emph{$F$-eigenvector} of $c$ if $v^Fc= v$ 
(note that the ``$F$-eigenvalue" is always one).
The set $E(k)$ of all $F$-eigenvectors in $V$ is a $k$-space of dimension 
$d$.
By Lang's theorem, the $k_{rs}$-span of $E(k)$ must be
equal to $V(k_{rs})$.
There is a variety $E$ defined over $k$ such that $E(k_t)$ is the 
$k_t$-span of $E(k)$ for every positive integer $t$.
Such a variety is called a \emph{$k$-form} of~$V$
\cite[Section~11.1]{Springer98Book}.

The following easy lemma is the key to our recursive approach.
\begin{lemma}
Let $G$ be a connected linear algebraic group defined over $k$ and
let $V$ be a $G$-module defined over $k$ with kernel $K\le G$.
Let $c$ be an element of $G$.
Suppose that $E(k)$ is the set of  $F$-eigenvectors of $c$ in $V$.
Then $a\in G$ satisfies $c\in a^{-F}Ka$ if, and only if,
$V(k)a=E(k)$.
\end{lemma}
\begin{proof}
If $a^{-F}za=c$ for $z\in K$, then, for all $v\in V(k)$,
$va=vza=va^Fc=(va)^Fc$ and so $va\in E(k)$.
Conversely, if $V(k)a=E(k)$, then, for all $v\in V(k)$,
$va=(va)^Fc=va^Fc$ and so $a^Fca^{-1}\in K$.
\end{proof}

We call an element $a\in G(k_{rs})$ such that
$V(k)a=E(k)$ a \emph{transformer} in $G$ for the $k$-form $E$.
Our approach to solving Lang's Theorem is outlined in 
Algorithm~\ref{A-lang}.
\begin{figure}
\begin{tabbing}
\quad\=\quad\=\quad\=\quad\=\kill
$\text{\sc Lang} := \text{\bf function}(G,c,s)\qquad$ [$c\in G(k_r)$, $s$ a multiple of the
order of $c^{F^{r-1}}\cdots c^Fc$]\\
\> \tbf{construct} a module $V$ for $G$\\
\> \tbf{let} $E(k)=\text{\sc $F$-Eigenspace}(c,V,s)$\\
\> \tbf{find} a transformer $a\in G(k_{rs})$ for $E$\\
\> \tbf{if} $V$ is faithful \tbf{then}\\
\> \> \tbf{return} $a$\\
\> \tbf{else}\\
\> \> \tbf{construct} a proper connected subgroup $H$ of $G$ containing the kernel of $V$\\
\> \> \tbf{let} $b=\text{\sc Lang}(H,a^Fca^{-1},s)$\\
\> \> \tbf{return} $ba$\\
\> \tbf{end if}\\
\textbf{end function}
\end{tabbing}
\caption{Algorithm outline for Lang's Theorem}
\label{A-lang}
\end{figure}
Note that $s$ is taken to be the order of $c^{F^{r-1}}\cdots c^F c$ in the 
top-level function call.
It is not necessary to recompute $s$ for the
recursive calls since a multiple of the element order works just as well.

Suppose that $G$ is split reductive and let $T_0$ be a
$k$-split maximal torus of $G$.
Using the methods of \cite{CohenMurrayTaylor04},
we can construct a module $V$ which is \emph{projectively faithful}, that is, 
the kernel $K$ is contained in the centre of $G$.
We can now take $H=T_0$ in Algorithm~\ref{A-lang}, since $Z(G)$ is contained in
every maximal torus of $G$.
Since a split torus has an easily constructed faithful module,
there is at most one recursive call for reductive groups.
The same algorithm could, in principle, be used for a nonreductive connected
group $G$: construct a simple quotient $G/N$, take $V$ to be the $G$-module
induced by a projectively faithful module for $G/N$, and take $H$ to be the
preimage in $G$ of the maximal torus in $G/N$.
However, finding the normal subgroup $N$ and constructing the quotient $G/N$ are
nontrivial problems which lie beyond the scope of this paper.

Algorithms for finding transformers are discussed in the next section.
We now give two algorithms for computing the $F$-eigenspace.
The most straightforward method is given in Algorithm~\ref{A-evectdet}.
\begin{figure}
\begin{tabbing}
\quad\=\quad\=\quad\=\quad\=\kill
$\text{\sc $F$-eigenspace} := \text{\bf function}(c,V,s)\qquad$ [$c\in G(k_r)$, $s$ the
order of $c^{F^{r-1}}\cdots c^Fc$]\\
\>  \tbf{let} $S$ be the $k$-matrix of $F$ acting of $k_{rs}$\\
\>  \tbf{let} $C$ be the $k$-matrix of $c$ acting on $V(k_{rs})=k_{rs}^{\;d}$\\
\>  \tbf{return} the fixed point space of $S^{\ds d}C$\\
\textbf{end function}
\end{tabbing}
\caption{Deterministic method for computing $F$-eigenvalues}
\label{A-evectdet}
\end{figure}
The key is to consider $k_{rs}$ as a $k$-space of dimension $rs$ and to consider
$V(k_{rs})= k_{rs}^{\;d}$ as a $k$-space of dimension $drs$.
The solution is then found by linear algebra over $k$.
Computing $S$ takes time $O(r^2s^2\log^2(q))$, where the second factor of $\log(q)$
is for computing $q$th powers.
Finding $C$ and the fixed space takes time $O(d^3r^3s^3\log(q))$.
So the overall time is $O(d^3r^3s^3\log^2(q))$.

An alternative method, due to Glasby and Howlett \cite{GlasbyHowlett97},
is given in Algorithm~\ref{A-evectlv}.
\begin{figure}
\begin{tabbing}
\quad\=\quad\=\quad\=\quad\=\kill
$\text{\sc $F$-eigenspace} := \text{\bf function}(c,V,s)\qquad$ [$c\in G(k_r)$, $s$  the
order of $c^{F^{r-1}}\cdots c^Fc$]\\
\>  \tbf{repeat}\\
\>  \>  \tbf{let} $x$ be a random $d\times d$ matrix over $k_{rs}$\\
\>  \>  \tbf{let} $a=x + x^Fc + x^{F^2}c^Fc + \cdots + x^{F^{rs-1}}c^{F^{rs-2}}\cdots c^Fc$\\
\>  \tbf{until} $a$ is invertible\\
\>  \tbf{return} $V(k)a^{-1}$\\
\textbf{end function}
\end{tabbing}
\caption{Las Vegas method for computing $F$-eigenvalues}
\label{A-evectlv}
\end{figure}
It takes time $O(d^2r^2s^2\log^2(q))$ to apply $F$ to a $d\times d$ matrix over $k_{rs}$,
so computing $a$ takes time $O(d^3r^2s^2\log^2(q))$.
Each randomly chosen $x$ has a probability of at least $1/4$ of yielding an
invertible element $a$.
Since this probability is bounded away from zero as $q$, $r$, $s$, and $d$
become large, the algorithm is Las Vegas. 
Note that we have an algorithm for Lang's theorem in $\gl{}{V}$
if the function returns $a$ instead of $V(k)a^{-1}$.

We now have:
\begin{theorem}\label{T-Etime}
Let $G$ be a connected linear algebraic group defined over $k$ and
let $V$ be a $G$-module defined over $k$ with dimension $d$.
Let $c$ be an element of $G$ with minimum field degree $r$ and
let $s$ be the order of $c^{F^{r-1}}\cdots c^Fc$.
Then we can compute a basis for the $k$-space $E(k)$ of $F$-eigenvectors of $c$ in 
deterministic time $O(d^3r^3s^3\log^2(q))$ 
or Las Vegas time $O(d^3r^2s^2\log^2(q))$.
\end{theorem}

\section{Finding transformers}\label{S-transf}
Let $G$ be a $k$-split connected reductive linear algebraic group defined over
$k$, let $c$ be in $G(k_r)$, and let $s$ be the order of 
$c^{F^{r-1}}\cdots c^Fc$.
Let $T_0$ be the standard $k$-split maximal torus of $G$ determined by the
Steinberg presentation.
Let $V$ be a projectively faithful $G$-module 
and compute $E$, the $k$-form of $F$-eigenvectors of $c$. 
In this section, we show how to find a transformer $a\in G(k_{rs})$
such that $E(k)a=V(k)$.
First we consider two special cases: split tori and classical groups.
Then we give an algorithm for an arbitrary
$k$-split connected reductive group.
The key is to consider $k$-bases with some additional structure that ensures that
$G$ is transitive on all such bases (or $G_{\rm ad}$ is transitive in
Subsection~\ref{SS-chevtransf}).  

\subsection{Split tori and isogeny}
A \emph{$k$-split torus} $T$ of dimension $n$ is just the direct product of
$n$ copies of $\kbar^\times$ with the Frobenius endomorphism taking the $q$th
power of each component. 
The standard module $V$ is just $\kbar^n$ with the componentwise action.
Suppose $c=(c_1,\dots,c_n)\in T(k_r)$ and $E$ is the variety of $F$-eigenvectors
of $c$ in $V$.
Splitting $T$ into $n$ copies of $\kbar^\times$
and using Theorem~\ref{T-Etime}, we can compute $E$ in Las Vegas time
$O(nr^2s^2\log^2(q))$.
Now $E$ has a basis of the form $a_1e_1,\dots,a_ne_n$
where each $a_i\in k_{rs}^{\;\times}$ and $e_i$ is the $i$th standard basis
vector in $V$.
Now $(a_1,\dots,a_n)\in T(k_{rs})$ is a transformer for $E$.
Hence we have proved:
\begin{proposition}\label{P-torustime}
Let $T$ be a $k$-split torus of dimension $n$.
Let $c$ be in $T(k_r)$, and suppose we are given $s$, the order of
$c^{F^{r-1}}\cdots c^Fc$.
Then we can find an element $a$ in $T(k_{rs})$ such that $c=a^{-F}a$ in
Las Vegas time $O(nr^2s^2\log^2(q)).$
\end{proposition}

Consider two connected linear algebraic groups $G$ and $H$ defined over $k$.
Let $\iota$ be a homomorphism $G\to H$ defined over $k$ which is onto
with finite kernel $K$.
Such a map is called an \emph{isogeny}.
Now suppose $G$ and $H$ are reductive and described by a Steinberg presentation
with unipotent, Weyl, and toral generators as in \cite{CohenMurrayTaylor04}.
If there is an isogeny $\iota:G\to H$, then both groups have the same Cartan
type.
Furthermore, we can assume (after composing with an automorphism)
that $\iota$ leaves unipotent and Weyl generators
unchanged, and acts by a change of basis on the toral generators.
We denote the standard tori generated by the toral generators of $G$ (resp.\
$H$) by $T_0$ (resp.\ $U_0$).
Note that $K\le Z(G)\le T_0$.
An important invariant of $\iota$ is the exponent of the finite group 
$K$, which we denote by $m$.

For $g\in T_0(k_r)$, we have $\iota(g)^{F^r}=\iota(g^{F^r})=\iota(g)$, so
$\iota(g)\in U_0(k_r)$.
This image can be computed in time $O(n^3r\log(q))$ by linear algebra in
$T_0(k_r)$.

For $h\in U_0(k_r)$, we can find $g\in T_0$ such that $\iota(g)=h$.
Then $\iota(g^{-F^r}g)=h^{-F^r}h=1$, ie, $g^{-F^r}g\in K$.
Hence $(g^{-F^r}g)^m=1$ and so $g^m\in T_0(k_r)$.
Using the fact that $T_0$ is a direct sum of copies of $\kbar^\times$,
such a $g$ must be in $T_0(k_{rm})$.
This preimage can be computed in time $O(n^3rm\log(q))$ by linear algebra in
$T_0(k_{rm})$.

\begin{proposition}\label{P-isogeny}
Let $G$ and $H$ be $k$-split connected reductive linear algebraic groups
defined over $k$ with reductive rank $n$.
Suppose we have an isogeny $\iota:G\to H$ defined over $k$ and
that $m$ is the exponent of the kernel of $\iota$.
For $c$ in $G$ or $H$, let $s(c)$ denote the order of $c^{F^{r-1}}\cdots c^F c$, where
$r$ is the minimal field degree of $c$.
\begin{enumerate}
\item\label{P-isogeny-GfromH}
  Lang's theorem can be solved for $c\in G(k_r)$ in time
  $O(n^3r^2s(c)^2m^2\log^2(q))$ plus the time needed to 
  solve it for some $c'\in H(k_r)$ with $s(c')\le s(c)$.
\item\label{P-isogeny-HfromG}
  Lang's theorem can be solved for $c\in H(k_r)$ in time 
  $O(n^3rs(c)m^2\log(q))$ plus the time needed to solve it for some
  $c'\in G(k_{rm})$ with $s(c')\le ms(c)$.
\end{enumerate}
\end{proposition}
\begin{proof}
If $c\in G(k_r)$, then $c'=\iota(c)$ can be found in time $O(n^3r\log(q))$.
Clearly $s(c')\le s(c)$.  Now we can find $a'\in H_{rs(c)}$ such that
$a'^{-F}a'=c'$.
Let $a\in G_{rsm}$ be a preimage of $a'$ computed in time $O(n^3rm\log(q))$.
Consider $a^Fca^{-1}\in K(k_{rs(c)m})\le T_0(k_{rs(c)m})$.
Now
$$
  (a^Fca^{-1})^{F^{rs(c)m-1}}\cdots(a^Fca^{-1})^F(a^Fca^{-1})
    = a^{F^{rs(c)m}}(c^{F^{rs(c)m-1}}\cdots c^F c)a^{-1}=1,
$$
So by Proposition~\ref{P-torustime},
we can find $b\in T_0(k_{rs(c)m})$ such that
$a^Fca^{-1}=b^{-F}b$ in Las Vegas time $O(nr^2s(c)^2m^2\log^2(q))$.
Now $(ab)^{-F}ab=c$ and Part~(1) follows.

If $c\in H(k_r)$, we can find an element $c'\in G(k_{rm})$ such that 
$\iota(c')=c$ in time $O(n^3rm\log(q))$.
Since
$(c'^{F^{r-1}}\cdots c'^Fc')^{s(c)}\in K$, we get $s(c')\le ms(c)$.
We can now find $a'\in H(k_{rs(c)m^2})$ such that $c'=a'^{-F}a'$.
Then $a=\iota(a')$ can be computed in time $O(n^3rs(c)m^2\log(q))$
and $a^{-F}a=c$ and Part~(2) is proved.
\end{proof}

\subsection{Classical groups}\label{SS-classtime}
We now show how to find transformers for the classical groups, using the
standard representations.
Throughout this subsection we take $V=\kbar^d$ and 
$B_0$ to be the standard basis $e_1,\dots,e_d$ of $V(k)$.

The simplest case is $G=\gl{d}{\kbar}$.
Let $B$ be a $k$-basis of $E(k)$.
Let $a$ be the matrix whose rows are the elements of $B$.
Then $B_0a=B$, and so $a$ is a transformer for $E$.

Now suppose $G=\sl{d}\kbar$.
Given a basis $B$ of $V$, define its \emph{volume}, 
denoted $\vol(B)$, to be the determinant of the matrix whose rows are the
elements of $B$. 
Then $B_0$ has volume one and $G$ is transitive on all
bases of volume one.
Now suppose $B$ is a basis of $E(k)$, the set of $F$-eigenvectors of $c\in G$.
Then $B^Fc=B$, and so
$$
  \vol(B)^F = \vol(B^F) = \vol(Bc^{-1}) = \vol(B) \det(c)^{-1} = \vol(B),
$$
and $\vol(B)\in k$.
We can now construct a basis $B'$ of $E(k)$ with volume one by
dividing the first element of $B$ by the scalar $\vol(B)$.
So the matrix that takes $B_0$ to $B'$ is a transformer in $G$.

Now suppose that $q$ is odd and $M$ is a nondegenerate orthogonal or symplectic
form on $V$ written $(u,v)$ for $u,v\in V$.
Further suppose that $M$ is defined over $k$.
Then the invariant group 
$$
  G = \{ x \in \GL_d(\kbar) \mid (ux,vx)=(u,v) \}
$$
is a reductive linear algebraic group defined over $k$.
Note that $G$ is not necessarily split or connected however.
Define the $m\times m$ matrix
$$
  A_m = \left(\begin{matrix} 0& & 1 \\ &\addots& \\1& & 0\end{matrix}\right)
$$
and let $\del$ be a fixed nonsquare in $k$.
Then the form $M$ has precisely one of the following Gram matrices $M_B$
with respect to some basis $B$:
\begin{itemize}
\item If $M$ is orthogonal and $d=2\ell+1$, then
  $$
    M_B=A_d \quad\text{or}\quad
    \left(\begin{matrix} & & A_\ell \\
                       & \del& \\
                       A_\ell&& \end{matrix}\right).
  $$
\item If $M$ is orthogonal and $d=2\ell$, then
  $$
    M_B=A_d \quad\text{or}\quad
    \left(\begin{matrix} & & & A_{\ell-1} \\
                       &1& & \\
                       & &-\del& \\
                       A_{\ell-1}& & \end{matrix}\right).
  $$
\item If $M$ is symplectic, then $d=2\ell$ and
  $$
    M_B=\left(\begin{matrix} & A_\ell \\ & \\-A_\ell & \end{matrix}\right).
  $$
\end{itemize}
A \emph{normal basis} for $M$ is a basis of $V$ such that $M_B$ is one of these
matrices. 

Given a nondegenerate symplectic or orthogonal form $M$ on the $k$-space $U$, 
Algorithm~\ref{A-sympl} constructs a
normal basis for $U$.
\begin{figure}
\begin{tabbing}
\quad\=\quad\=\quad\=\quad\=\kill
$\text{\sc NormalBasis} := \text{\bf function}(U)$\\
\> \tbf{let} $u$ be a nonzero element of $U$\\
\> \tbf{if} $\dim(U)=1$ \tbf{then}\\
\> \> \tbf{find} $a\in k$ such that $a^2=(u,u)$ or $a^2\del=(u,u)$\\
\> \> \tbf{return} $u/a$\\
\> \tbf{end if}\\
\> \tbf{let} $v$ be a nonzero element of $u^\perp\setminus ku$\\
\> \tbf{if} $\dim(U)=2$ \tbf{then}\\
\> \> \tbf{if} $(u,u)\in -\del(v,v) k^{\times2}$ \tbf{then} $\qquad$[$U$ anisotropic]\\ 
\> \> \> \tbf{find} $a,b\in k$ such that $(u,u)a^2+(v,v)b^2=1$\\
\> \> \> \tbf{find} $c\in k$ such that $((u,u)a^2-(v,v)b^2)c^2=-\del$\\
\> \> \> \tbf{return} $au+bv, c(au-bv)$\\
\> \> \tbf{else} $\qquad$[$U$ isotropic]\\
\> \> \> \tbf{find} $a,b\in k$ such that $(u,u)a^2+(v,v)b^2=0$\\
\> \> \> \tbf{let} $w$ be a nonzero element of $(au+bv)^\perp$\\
\> \> \> \tbf{return} $au+bv$, $w$\\
\> \> \tbf{end if}\\
\> \tbf{end if}\\
\> \tbf{let} $w$ be a nonzero vector in $(ku+kv)^\perp\setminus(ku+kv)$\\
\> \tbf{find} $a,b,c\in k$ such that $(u,u)a^2+(v,v)b^2+(w,w)c^2=0$\\
\> \tbf{let} $x$ be a nonzero element of $(au+bv+cw)^\perp$\\
\> \tbf{return} $au+bv+cw, \text{\sc NormalBasis}(\{u,v\}^\perp), x$\\
\tbf{end function}
\end{tabbing}
\caption{Finding a normal basis for a space with a bilinear form}
\label{A-sympl}
\end{figure}
The quadratic equations involved always have solutions by the standard
classification theory of bilinear forms over finite fields (see 
\cite{Grove02Book} for more details).  Each of these equations
can be solved by standard techniques in time $\log^2(q)$.
Note that this construction is rational (that is, it does not use extensions of 
$k$) and takes time $O(d^3\log(q)+d\log^2(q))$.

If the form $M$ is symplectic, we are done: our transformer
is simply the matrix taking a normal basis of $V(k)$ to a normal basis of $E(k)$.

If $M$ is orthogonal, the two normal bases may have different Gram matrices,
in which case the equation in Lang's Theorem has no solution.
This is to be expected, since $G$ is not connected in this case.
If we take $G=\SO(V,M)$,
then this problem can be avoided.
Without loss of generality, the standard basis $B_0$ is normal.
Suppose that $B$ is a normal basis of $E(k)$.
As with the special linear group, $\vol(B)$ is in $k$.
Also
$$
  \det(M_{B}) = \det(BM_{B_0}B^T)=\vol(B)^2\det(M_{B_0}).
$$
But the two choices given above for the Gram matrix have determinants in
different cosets of $k^{\times2}$, hence $M_{B_0}=M_B$.
It now remains to ensure that $\vol(B)$ is one.
Now $\vol(B)^2=\det(M_{B_0})/\det(M_B)=1$, so suppose $\vol(B)=-1$.
If $M_B=A_d$, then exchanging the first and last vectors in $B$ results in a new
normal basis with volume one.
Otherwise, negating the $(\ell+1)$st vector in $B$ has the same effect.

Similar methods work for quadratic forms in characteristic two.  

Now suppose we have a split simple classical group $G$ of (reductive and
semisimple) rank $n$.
Then $G$ is isogenous to one of the groups considered above, with $d=O(n)$.

If $G$ has type $\CA_n$, then there is an isogeny map 
$\iota:\sl{n}{\kbar}\to G$ with $m\le n+1$.
By Proposition~\ref{P-isogeny}(\ref{P-isogeny-HfromG}),
we can solve Lang's Theorem in $G$ in Las Vegas time $O(n^5r^2s^2\log^2(q))$.

If $G$ is not of type $\CA_\ell$, then there is a series of at most 2
isogeny maps connecting $G$ with one of the groups considered above.
For each of these maps, $m$ is at most 4.
By Proposition~\ref{P-isogeny}, we can solve Lang's Theorem in Las Vegas time
$O(n^3r^2s^2\log^2(q))$.

We have now proved Theorem~\ref{T-classtime}.
For groups of Cartan type $\CG_2$ and $\CF_4$, a
similar result can probably be obtained by exploiting the structure of
composition and Jordan algebras, respectively \cite{SpringerVeldkamp00Book}.

\subsection{Adjoint representation}\label{SS-chevtransf}
Now consider an arbitrary $k$-split connected reductive linear algebraic group $G$, 
with reductive rank $n$ and semisimple rank $\ell$. 
Then $G$ has a root datum $(X,\Phi,Y,\Phi^\star)$ with respect to a $k$-split
maximal torus $T_0$.
Here $X$ and $Y$ are free $\Z$-modules of dimension $n$ with a bilinear
pairing $\langle\circ,\circ\rangle:X\times Y\to \Z$ putting them in duality.
We fix dual bases $e_1,\dots,e_n$ for $X$ and $f_1,\dots,f_n$ for $Y$.
The roots $\Phi$ are a finite subset of $X$ and the coroots 
$\Phi^\star$ are a finite subset of $Y$.
There is a one-to-one correspondence $\star:\Phi\to\Phi^\star$ 
such that $\langle\al,\al^\star\rangle=2$ for every $\al\in\Phi$  .
For more details see \cite{CohenMurrayTaylor04}.

The Lie algebra $L=L(G)$ is a $G$-module defined over $k$.
This is called the \emph{adjoint representation} of $G$ and
it is projectively faithful.
Now $L(k)$ has basis elements $e_\al$ for $\al\in\Phi$ and $h_i\in L(T_0)$ for
$i=1,\dots,n$ with structure constants:
\begin{align}
  [h_i,h_j]   &= 0,\label{E-hh}\\
  [e_\alpha,h_i] &= \langle\alpha,f_i\rangle\,e_\alpha,\label{E-eh}\\
  [e_{-\alpha},e_\alpha] &= \sum_{i=1}^n \langle e_i,\al^\star\rangle
  h_i,\label{E-em}\\
  [e_\alpha,e_\beta]   &= \begin{cases}
    N_{\alpha\beta}\,e_{\alpha+\beta}&\text{for
    }\alpha+\beta\in\Phi,\\
    0&\text{for $\al+\be\notin\Phi$, $\be\ne-\al$,}
  \end{cases}\label{E-ee}
\end{align}
where the integral constants $N_{\al\be}$ are defined in
\cite{Carter72Book}.
Such a basis is called a \emph{Chevalley basis}.

Choose simple roots $\al_1,\dots,\al_\ell$, and fix a linear ordering $<$
on $\Phi^+$ which is compatible with height, 
ie, $\height(\al)<\height(\be)$ implies that $\al<\be$.
Given a nonsimple positive root $\xi$, take the positive roots $\al,\be$ such that $\xi=\al+\be$
and $\al$ is as small as possible with respect to the ordering on $\Phi^+$.
We call $(\al,\be)$ the \emph{extraspecial pair} of $\xi$.
We can choose a Chevalley basis of $L$ so that the integers $N_{\al\be}$ are
positive on extraspecial pairs by \cite{Carter72Book}.
We call such a basis a \emph{standard Chevalley basis}.
Note that, as with the normal bases in Subsection~\ref{SS-classtime},
the problem of finding a standard Chevalley basis is rational over $k$.

The linear map $a$ taking the standard Chevalley basis of $L(k)$ to
a standard Chevalley basis of $E(k)$ must be an automorphism of $L(k_{rs})$.
We now need to find a transformer in $G$.
Let $G_{\rm ad}$ be the adjoint group with the same Cartan type as $G$ and
let $\Gamma$ be the automorphism group of the Dynkin diagram of $G$.
For each element of $\Gamma$, fix a graph automorphism normalising $T_0$ and the
Borel subgroup determined by $\Phi^+$, as in \cite{Carter72Book}.
Take $Z$ to be a complement to
$Z(L)\cap[L,L]$ in $Z(L)$; by construction, our graph automorphisms fix $Z$
pointwise.
If the characteristic of $k$ is greater than $3$, then it follows from
\cite{Hogeweij82} that
$$
  \Aut(L) = \Aut(Z)\times(\Gamma\ltimes G_{\rm ad}) .
$$
We can compute a decomposition $a=z\ga b$ with $z\in\Aut(Z)$, $\ga$ a graph
automorphism, and
$b\in G_{\rm ad}(k_{rs})$  in time $O(d^3rs\log(q))$
using a slight modification of Algorithm~6 of \cite{CohenMurrayTaylor04}.
Since $L(k)z\ga=L(k)$, the element $b$ is a transformer in $G_{\text{ad}}$.
It is easily checked on a case-by-case basis that the number of roots of $G$ is
$O(\ell^2)$ and so the dimension of $L$ is $O(n+\ell^2)$.
Hence Lang's Theorem can be solved for $c\in G_{\rm ad}(k_r)$ in time
$O((n+\ell^2)^3r^2s^2\log^2(q))$, once we have a standard Chevalley basis for
$E(k)$.

Suppose now that $G$ is simple, so that $\ell=n$.
Then there is an isogeny map $G\to G_{\rm ad}$ with $m$ at most $n+1$. 
We can now apply
Proposition~\ref{P-isogeny}(\ref{P-isogeny-HfromG}) and obtain:
\begin{proposition}\label{P-ad}
Suppose that $k$ has characteristic greater than $3$.
Let $G$ be a $k$-split connected simple linear algebraic group
and let $L$ be the Lie algebra of $G$.
Let $c$ be an element of $G(k_r)$ and suppose we are given $s$, the order
of $c^{F^{r-1}}\cdots c^Fc$.
Let $E$ be the variety of $F$-eigenvectors of $c$ in $L$.
We can find $a\in G(k_{rs})$ in Las Vegas time
$O(n^8r^2s^2\log^2(q))$ plus the time needed to find a standard Chevalley basis
of $E(k)$.
\end{proposition}

We give an algorithm for finding a standard Chevalley basis in the next section.
The timing of this algorithm is analysed in Section~\ref{S-prob}, leading to a
proof of Theorem~\ref{T-time}.

\section{Computing a standard Chevalley basis}\label{S-csa}
We now give an algorithm for constructing a Chevalley basis
of the Lie algebra $L$ of a $k$-split connected reductive group $G$.
Recall that $L$ is a $p$-Lie algebra \cite[Section~V.7]{Jacobson62Book}.
The first and hardest step is finding a $k$-split maximal toral $p$-subalgebra.
This is similar to the algorithm of \cite{GraafIvanyosRonyai96} for finding a
Cartan subalgebra, but ensuring that the subalgebra is $k$-split
makes things considerably more complex.
Once we have a split maximal toral $p$-subalgebra, a Chevalley basis can be
constructed using \cite[Section~4.2]{Carter72Book}.

Our algorithm only works for fields of characteristic $p>3$.
Whenever possible we state results for characteristics $2$ and $3$, in the hope
that the gaps in our argument for small $p$ can be filled later.

We assume that $L$ is given as a structure constant algebra, but we
frequently compute in the adjoint representation.
Throughout this section $n$ denotes the reductive rank of $G$, $\ell$ denotes
the semisimple rank of $G$, and $d$ denotes the dimension of the Lie algebra
$L$.
Recall that our Steinberg presentation of $G$ determines a
$k$-split maximal torus $T_0$.

\subsection{Toral subalgebras}\label{SS-toral}
A Lie algebra $L$ over a field of characteristic $p$ is called a 
\emph{$p$-Lie algebra} if it is equipped with a map $p:L\to L$ satisfying the
axioms 
\begin{align}
  (x+y)^p &= x^p+y^p + \sum_{i=1}^{p-1}s_i(x,y),\label{E-sum}\\
  (ax)^p  &= a^px^p,\label{E-sprod}\\
  [xy^p]  &= x(\ad y)^p\label{E-ad}
\end{align}
where $x,y\in L$, $a\in\kbar$, $s_i$ is defined in
\cite[Section~V.7]{Jacobson62Book}, 
and $a^p$ and $(\ad y)^p$ are the usual $p\/$th powers.

Given values of the $p$-map on a basis of $L$, we can compute the values
on an arbitrary element using Equations~(\ref{E-sum}) and~(\ref{E-sprod}).
But $s_{p-1}$ involves Lie products of length $p$, so the time taken for this 
computation is not polynomial in $\log(p)$.
Given $x\in L$, we can use (\ref{E-ad}) to compute the coset $x^p+Z(L)$ in
time $O(\ell^6\log(q)\log(p))$, since $\dim(L/Z(L))$ is $O(\ell^2)$.
We also define the \emph{$q$-map} by applying the $p$-map $e$ times, where $q=p^e$;
values of this map modulo $Z(L)$ can be computed in time 
$O(\ell^6\log^2(q))$.

We say that $x\in L$ is \emph{semisimple} if it is contained in the
$p$-subalgebra generated by $x^p$.
A \emph{toral subalgebra} of $L$ is a subalgebra defined over $k$ consisting
entirely of semisimple elements.
Note that a toral subalgebra need not be a $p$-subalgebra.
However every subalgebra $H$ of $L$ is contained in a minimal $p$-subalgebra
called the \emph{$p$-closure} of $H$ in $L$.
The $p$-closure of a toral subalgebra is toral, and so a maximal
toral subalgebra is automatically a $p$-subalgebra.
An $n$-dimensional toral $p$-subalgebra $H$ is \emph{$k$-split} if
$H(k)$ is isomorphic, as a $p$-Lie algebra, to the vector space $k^n$
with trivial Lie product and the $p$-map acting componentwise.

If $L$ is the Lie algebra of a $k$-split connected reductive group $G$, then
the values of the $p$-map on a Chevalley basis are
$$
  {h_i}^p=h_i \quad\text{and}\quad {e_\al}^p=0,$$
provided that $p>3$.
Clearly 
$
  H_0:= L(T_0)=\langle h_1,\dots,h_n\rangle
$
is a $k$-split toral subalgebra.

We say that the Lie algebra $L$ is \emph{$k$-split} if it contains a
maximal toral subalgebra which is $k$-split.
The following theorem collects together the properties of toral
subalgebras which we need. 
\begin{theorem}\label{T-toral}
  Let $L$ be the $p$-Lie algebra of a $k$-split connected reductive group $G$.
  \begin{enumerate}
    \item\label{T-toral-split} $L$ is $k$-split with split maximal toral subalgebra
      $H_0$.
    \item\label{T-toral-centre} The centre of $L$ is a $k$-split toral subalgebra when
    $p>2$.
    \item\label{T-toral-abelian} Every toral subalgebra of $L$ is abelian.
    \item\label{T-toral-torus} Every ($k$-split) maximal toral subalgebra of $L$ 
      is the Lie algebra of a ($k$-split) maximal torus of $G$ (when $p>2$).
    \item\label{T-toral-conj} The maximal toral subalgebras of $L$ are
      $G$-conjugate.
    \item\label{T-toral-splitconj}  The $k$-split maximal
      toral subalgebras of $L$ are $G(k)$-conjugate when $p>2$.
  \end{enumerate}
\end{theorem}
\begin{proof}
In Part~(\ref{T-toral-split}) it only remains to prove maximality, which follows from 
\cite[Proposition~13.2]{Humphreys67Book}.
Part~(\ref{T-toral-abelian}) is given in \cite[Lemma 8.1]{Humphreys78Book}
for characteristic zero, but the same proof works for positive characteristic.
Part~(\ref{T-toral-conj}) is Corollary~13.5 of \cite{Humphreys67Book}. 

We now prove Part~(\ref{T-toral-centre}).
Let $\{e_\al,h_i\}$ be a Chevalley basis with respect to $H_0$.
Suppose that $z\in Z(L)$ and write 
$$
  z = \sum_{i=1}^n t_i h_i + \sum_{\al\in\Phi} a_\al e_\al.
$$
Let $h_\al=\sum_{i=1}^n\langle e_i,\al^\star\rangle h_ i$, then the coefficient of $e_\al$ in
$[z,h_\al]$ is $2a_\al$.
Since $[z,h_\al]=0$ and $p>2$, we get $a_\al=0$.
Hence $z$ is in $H_0=\langle h_1,\dots,h_n\rangle$.
Since $H_0$ is a split toral subalgebra, $Z(L)$ is also.
(The idea for this proof is from \cite[Lemma~6.10]{Hogeweij82}.) 

Every maximal toral subalgebra of $L$ is the Lie algebra of a maximal torus of
$G$ by \cite[Proposition~13.2]{Humphreys67Book}.
For $p>2$, split tori correspond to split toral subalgebras by
\cite[Theorem~9]{Seligman67}.
Hence Part~(\ref{T-toral-torus}) is proved.

From now on we assume $p>2$.
By \cite[Proposition~13.6]{Humphreys67Book}, $T\goesto L(T)$ gives a
one-to-one correspondence between maximal tori of $G$ and maximal toral
subalgebras of $L$.
Once again, split tori correspond to split toral $p$-algebras.
Part~(\ref{T-toral-splitconj}) now follows from the corresponding result for
tori.
\end{proof}

\begin{corollary}
Given a subalgebra $H$ of $L$ defined over $k$, we can determine if $H$ is $k$-split maximal
toral in time $O(\ell^7\log^2(q))$.
\end{corollary}
\begin{proof}
First check that $H$ is abelian of dimension $n$, and compute $H/Z$.
Note that $H/Z$ has dimension at most $\ell$.
By Theorem~\ref{T-toral}(\ref{T-toral-centre}), it suffices to determine if
$H/Z(L)$ is a split toral algebra.
This is done by testing whether
$b^q+Z(L)=b+Z(L)$ where $b+Z(L)$ runs over a basis of $L/Z(L)$.
As we argued at the beginning of this section, 
this takes time $O(\ell^6\log^2(q))$ for each basis element.
\end{proof}

Since semisimple elements are common in $L(k)$ (see Section~\ref{S-prob}) and
the centraliser of such an element is reductive of rank $n$, we can find a 
maximal toral subalgebra by
Algorithm~\ref{A-toral}.
\begin{figure}
\begin{tabbing}
\quad\=\quad\=\quad\=\quad\=\kill
$\text{\sc MaximalToralSubalgebra}:=\text{\bf function}(L)$\\
\>  \textbf{repeat} take $x$ random in $L$ \textbf{until} $x$ is semisimple\\
\>  \textbf{let} $M=C_L(x)$\\
\>  \textbf{if} $M$ is abelian \textbf{then}\\
\>  \>  \textbf{return} $M$\\
\>  \textbf{else}\\
\>  \>  \textbf{return} $\text{\sc MaximalToralSubalgebra}(M)$\\
\>  \textbf{end if}\\
\textbf{end function}
\end{tabbing}
\caption{Finding a maximal toral subalgebra in $L$}
\label{A-toral}
\end{figure}

The basic idea of our algorithm is to randomly select a series of increasingly split
maximal toral subalgebras.
We now assign a conjugacy class of $W$ to every maximal toral subalgebra $H$,
which measures how split $H$ is.
See \cite{Lehrer92} for a more detailed version of this construction.
There exists $g\in\offd{G}{\kbar}$ such that $H={H_0}^g$, by
Theorem~\ref{T-toral}(\ref{T-toral-conj}).
Note that ${H_0}^F=H_0$ and $H^F=H$, since both are defined over $k$.
Now
$$
  {H_0}^{g^Fg^{-1}} = (({H_0}^g)^F)^{g^{-1}} = (H^F)^{g^{-1}} = H^{g^{-1}} =H_0,
$$
so $g^Fg^{-1}$ is in $N_G(H_0)=N_G(T_0)$.
Let $w$ be the image of $g^Fg^{-1}$ under projection onto the
Weyl group $W=N_G(T_0)/T_0$.
The element $w$ is uniquely determined by $H$ up to conjugacy in $W$.

\subsection{Root decompositions of $L$}\label{SS-rtdec}
The \emph{root decomposition} of $L$ with respect to $H_0$ is
$$
  L = H_0\ds \bigds_{\al\in\Phi} L_{\al}
$$
where the \emph{root space}
$L_{\al}=\{ b \in L \mid [b,h]=\al(h)b \text{ for all $h\in H_0$}\}$
and each root $\al\in\Phi$ is a linear functional $H_0\to \kbar$ defined over
$k$.
This decomposition is defined over $k$ by \cite[Theorem~6]{Seligman67}.
If the characteristic of $k$ is greater than 3, every root space has
dimension one.

Let $H$ be a maximal toral subalgebra of $L$,
fix $g\in\offd{G}\kbar$ such that $H={H_0}^g$ and let $w$ be the image of
$g^Fg^{-1}$ in $W$.
For $\al\in\Phi$, define $\al^g:H\to \kbar$ by $\al^g(h)=\al(h^{g^{-1}})$.
Then the root
decomposition with respect to $H$ is
$$
  L =H \ds \bigds_{\al\in\Phi} {L_{\al^g}},
$$
where ${L_{\al^g}}=\{ b \in L \mid [b,h]=\al^g(h)b \text{ for all $h\in H$}\} =
{L_\al}^g$.
This decomposition is not defined over $k$ in general.

Fix a basis $h_1,\dots,h_n$ of $H$ and let $f=(f_1,\dots,f_n)$ be a sequence of
irreducible polynomials in $k[X]$ with $f_i(X)\ne X$ for at least one $i$. 
Define
$$
  L_f=\{ y\in L \mid yf_i(\ad(h_i))=0\text{ for $i=1,\dots,\ell$}\}.
$$
If $L_f\ne0$, we call $f$ a \emph{generalised root} and
$L_f$ a \emph{generalised root space}.
The \emph{generalised root decomposition} of $L$ with respect to $H$ is
$$
  L =H \ds \bigds_{f\in\sF} L_f,
$$
where $\sF=\sF(L,H)$ is the set of generalised roots of $L$ with respect to $H$.
This decomposition is defined over $k$.
The generalised roots are computed by Algorithm~\ref{A-rts}.
\begin{figure}
\begin{tabbing}
\quad\=\quad\=\quad\=\quad\=\kill
$\text{\sc GeneralisedRoots}:=\text{\bf function}(L,H)$\\
\>  \tbf{let} $h_1,\dots,h_n$ be a basis of $H$\\
\>  \textbf{let} $\sF=\{ () \}$ and \textbf{define}  $L_{()} = L$\\
\>  \textbf{for} $i=1,\dots,n$ \textbf{do}\\
\>  \>  \textbf{let} $\sF'=\emptyset$\\
\>  \>  \tbf{for} $f\in\sF$ \tbf{do}\\
\>  \>  \>  \tbf{compute} $g$, the characteristic polynomial of $h_i$ on $L_f$\\
\>  \>  \>  \tbf{for} $f_i$ an irreducible factor of $g$ \tbf{do}\\
\>  \>  \>  \>  \tbf{add} $(f_1,\dots,f_{i-1},f_i)$ to $\sF'$ where $f=(f_1,\dots,f_{i-1})$\\
\>  \>  \>  \>  \tbf{define} $L_{(f_1,\dots,f_i)}=\{ x\in L_f \mid xf_i(\ad(h_i))=0\}$\\
\>  \>  \>  \tbf{end for}\\
\>  \>  \tbf{end for}\\
\>  \>  \tbf{let} $\sF=\sF'$\\
\>  \tbf{end for}\\
\>  \tbf{remove} $(X,\dots,X)$ from $\sF\qquad$ [since $L_{(X,\dots,X)}=H$]\\
\>  \textbf{return} $\sF$\\
\textbf{end function}
\end{tabbing}
\caption{Generalised roots}
\label{A-rts}
\end{figure}
Complete factorisation of a polynomial of degree $d$ over $k$ takes time at
most $O(d^3(\log(d)+\log(q))\log(q))$, as shown in \cite{VonZurGathenGerhard03}.
Factoring the characteristic polynomials $g$ is the dominant contribution to the running time
of this algorithm.
Since each $g$ has degree at most $d$, and the sum of the degrees of all the
$g$s is at most $nd$, the algorithm 
takes time $O(nd^3(\log(d)+\log(q))\log(q))$. 

In fact, we do not apply this algorithm directly to $L$, since we want our
time to depend on $\ell$ but not on $n$ (this is necessary for analysing the
recursion in Algorithm~\ref{A-mts}).
By Theorem~\ref{T-toral}(\ref{T-toral-centre}), the centre
$Z(L)$ is contained in $H$.
So we can construct a basis $h_1,\dots,h_n$ for $H(k)$ with
$\langle h_{1},\dots,h_m\rangle=Z(L)$ central for some $m\le n$.
Extend this to a basis $B$ of $L(k)$.
Let $\phi$ be the pullback
map $L/Z(L)\to L$ which takes $b+Z$ to $b$ for all $b\in B$.
Note that $\phi$ is a linear map, but need not be a Lie algebra map.
We compute  in $L/Z(L)$, since it has dimension $O(\ell^2)$ independent of
$n$, and the results are then transfered into $L$ via $\phi$.
However, $L/Z(L)$ need not be the Lie algebra of a group of Lie type, so
most of our theoretical results do not apply to this quotient.

Let $\sF$ be the set of generalised roots of $L/Z(L)$ with respect to $H/Z(L)$.
Given $f=(f_1,\dots,f_m)\in\sF$ define the sequence
$f'=(f_1,\dots,f_m,X,\dots,X)$ of length $n$.
It is now easy to see that $\phi((L/Z(L))_f)=L_{f'}$.
Hence the generalised root
decomposition of $L$ with respect to $H$ follows immediately once we have the
decomposition of $L/Z(L)$ with respect to $H/Z$.
Since the dimension of $L/Z(L)$ is $O(\ell^2)$, the decomposition
of $L/Z(L)$ can be computed in time $O(\ell^7(\log(\ell)+\log(q))\log(q))$.

Given a generalised root $f\in\sF(L,H)$, the subspace $L_f$ is a direct sum 
of components $L_{\al^g}$ of the root decomposition with respect to $H$.
So we can partition $\Phi$ into subsets $\Phi_f$ such that
$L_f=\bigds_{\al\in\Phi_f}L_{\al^g}$.
Define the \emph{degree} of $f$ to be the lowest common multiple of the degrees
of the $f_i$.
Given a generalised root $f$, we define 
$$
  f_-= ( (-1)^{\deg(f_1)}f(-X), \dots, (-1)^{\deg(f_n)}f(-X) ).
$$
Clearly $\Phi_{f_-}=-\Phi_f$.
Note that we can have $f=f_-$ when the degree of $f$ is greater than one.

We now prove some properties of the sets $\Phi_f$.
\begin{lemma}\label{L-orbit}
Let $f$ be a generalised root of $L=L(G)$ with respect to $H$.
\begin{enumerate}
\item\label{L-orbit-act} The action of $F$ on 
  $\{L_{\al^g}\mid\al\in\Phi_f\}$ is equivalent to the action of $w$ on
  $\Phi_f$. 
\item\label{L-orbit-orb} $\Phi_f$ is a union of orbits of $w$ on $\Phi$.
\item\label{L-orbit-fix} If $\deg(f)=1$, then $w$ acts trivially of $\Phi_f$.
  If in addition $q>3$, then $\Phi_f$ contains a single root.
\item\label{L-orbit-2} If $\deg(f)=2$ and $f=f_-$, then $w$ acts on $\Phi_f$ by
  negation.
\end{enumerate}
\end{lemma}
\begin{proof}
Write $g^Fg^{-1}=t\wdot$ for some $t\in T_0$.
Now
$$
  {L_{\al^g}}^F = {L_{\al}}^{gF} = {L_\al}^{F^{-1}gF} =
  {L_\al}^{g^F} = {L_\al}^{t\wdot g} = {L_{\al w}}^g = L_{(\al w)^g},
$$
and so Part~(\ref{L-orbit-act}) is proved.
Part~(\ref{L-orbit-orb}) follows since ${L_{f}}^F=L_{f}$.

Part~(\ref{L-orbit-fix}) holds because $L_f$ is a root space when $\deg(f)=1$. 

Suppose $\deg(f)=2$ and $f=f_-$.  Let $\al\in\Phi_f$.  
Then ${L_{\al^g}}^{F}=L_{(\al w)^g}$ and 
so $(\al w)^g(h_i)$ and $\al^g(h_i)$ are conjugate roots of $f_i$.
But if $\deg(f_i)=2$, then $f=f_-$ implies that the conjugate roots are negatives
of each other.
And if $\deg(f_i)=1$, then $f=f_-$ implies that the only root of $f_i$ is zero.
In either case $(\al w)^g(h_i)=-\al^g(h_i)$ and so $w$ acts by negation.
\end{proof}

\subsection{Fundamental subalgebras}
Now that we have the generalised root decomposition of $L$ with respect to $H$, we consider the
subalgebra $M_f$ generated by a generalised root space $L_f$.
Such subalgebras often turn out to be fundamental:
We define a \emph{(split) fundamental subgroup} of $G$ as a connected
reductive subgroup normalised by a (split) maximal torus.
A subalgebra $M$ of $L$ is \emph{(split) fundamental} if it is 
the Lie algebra of a (split) fundamental group.
This subgroup is denoted $G_M$.
Fundamental subalgebras clearly normalise a maximal toral subalgebra.

The most important properties of such algebras for our purposes are:
\begin{theorem}\label{T-fund}$\;$
Suppose that $k$ has characteristic greater than $3$.
Let $M$ be a fundamental subalgebra of $L$ normalised by the maximal toral
subalgebra $H$.
\begin{enumerate}
  \item\label{T-fund-torus}  $M\cap H$ is a maximal toral subalgebra of
    $M$.
  \item\label{T-fund-splitext}   If $H$ is a
  split maximal toral subalgebra, then $M$ is split fundamental. 
\end{enumerate}
\end{theorem}
\begin{proof}
We have $H\le C_{M+H}(H)\le C_L(H)=H$, so $M+H$ has root decomposition
\begin{align}\label{E-chev}
  M+H&= H\ds\bigds_{\be}M_\be,
\end{align}
where $\be$ runs over $\Phi(M+H,H)$, the set of roots of $M+H$ with respect to
$H$.
Suppose $m+h\in M_\be$ where $m\in M$ and $h\in H$.  Then, for all
$h'\in H$,
$[m,h']=[m+h,h']=\al(h')(m+h)$.  But $[m,h']\in M$, and so $h\in M$ and
$M_\be\le M$.
By intersecting (\ref{E-chev}) with $M$, using the fact that each $M_\be$ has
dimension one, we obtain the root decomposition
$$
  M= (H\cap M) \ds\bigds_{\be}M_\be,
$$
and so $H\cap M$ is a maximal toral subalgebra of $M$.

Finally if $M$ normalises $H$ and $H$ is split, then, by
Theorem~\ref{T-toral}(\ref{T-toral-torus}), $H=L(T)$ for some split maximal
torus $T$ of $G$.
Then $G_M$ normalises $T$ and Part~(\ref{T-fund-splitext}) is proved.
\end{proof}

Recall that the closure $\Psibar$ of $\Psi\subseteq\Phi$ is just the set of all
roots that can be  written as a sum of elements of $\Psi$.
Note that if $\Psibar$ is also closed under negation, it is a \emph{subsystem}.
If $\Psibar$ is a subsystem, we say $w\in W$ is \emph{inner} on $\Psibar$ if the action of $w$ on
$\Psibar$ is induced by an element of $W(\Psibar)$.
\begin{lemma}\label{L-Psi}
Suppose that $k$ has odd characteristic. 
Let $M$ be the subalgebra generated by $\sum_{\al\in\Psi}L_{\al^g}$,
where $\Psi$ is an orbit in $\Phi$ under the action of $w\in W$.
\begin{enumerate}
  \item\label{L-psi-fund} $M$ is fundamental or soluble.
  \item\label{L-psi-simpsol} If $M$ is fundamental, then $G_M$ is semisimple.
  \item\label{L-psi-neg} $M$ is fundamental if, and only if,
    $\Psibar$ is a subsystem.
  \item\label{L-psi-inner} If\/ $\Psibar$ is a subsystem and $w$ is
    inner on $\Psibar$, then $M$ is split fundamental.
\end{enumerate}
\end{lemma}
\begin{proof}
Since $[\sum_{\al\in\Psi}L_{\al^g},H]\le \sum_{\al\in\Psi}L_{\al^g}$, we
have $[M,H]\le M$ and so $M$ normalises $H$.
Since $[L_\al,L_\be]\le L_{\al+\be}$ (recalling that $L_0=H$), we have
$$
  M = (H\cap M) \ds \bigds_{\al\in\Psibar} L_{\al^g}.
$$

Let $\Psi=\Psi_1\cup\cdots\cup\Psi_m$ be the finest decomposition of
$\Psi$ into pairwise orthogonal subsets.
Then $\Psibar=\Psibar_1\cup\cdots\cup\Psibar_m$ is also an orthogonal
decomposition.
Clearly $w$ permutes the sets $\Psi_i$ and, since $w$ is transitive on $\Psi$,
it must be transitive on them.
Since $-\Psibar_1$ is never orthogonal to $\Psibar_1$, we either have
$-\Psibar_1=\Psibar_1$ or $-\Psibar_1$ is disjoint from
$\Psibar_1$.
By the transitivity of $w$, whichever of these cases holds for
$-\Psibar_1$, also holds for all $-\Psibar_i$.
In particular, $\Psibar$ is closed under negation if{f} $\Psibar_1$ is.
Let $\psi$ be the sum of all the elements of $\Psibar_1$.
Now $\Psibar_1$ is closed under negation if{f} $\psi=0$
(since $\psi=0$ implies
$-\al=\sum_{\be\in\Psibar_1,\be\ne\al}\be\in\Psibar_1$ for all $\al\in\Psibar_1$,
and the converse is trivial).
We define $M_i=H_i \ds \bigds_{\al\in\Psibar_i} L_{\al^g}$, where
$H_i$ is the subalgebra of $H$ generated by $[L_{\al^g},L_{-\al^g}]$ for
all $\al\in\Psibar_i$.
Note that $M=\bigds_iM_i$.

Suppose first that $\psi\ne0$.
Then the root subsystem generated by $\Psi_1$ is just
$\Psibar_1\cup-\Psibar_1$.
Since this root subsystem is irreducible, $\psi$ induces an ordering on it
which makes $\overline{\Psi_1}$ the set of positive roots.
Hence $M_1$ is just the Borel subalgebra of the Lie algebra of a simple group,
and so must be soluble.
The transitivity of $w$ on the sets $\Psibar_i$ implies that $M_i$ is soluble
for every $i$, and so $M=\bigds_iM_{i}$ is soluble.

If $\psi=0$, then $\Psibar_1$ is an irreducible root subsystem
and so $M_1$ is fundamental with $G_{M_1}$  a simple group.
Hence $M=\bigds_iM_{i}$ is fundamental with $G_M$ a semisimple group.
Parts~(\ref{L-psi-fund}), (\ref{L-psi-simpsol}) and~(\ref{L-psi-neg}) are now proved.

Now suppose that $\Psibar$ is closed under negation and $w$ is inner on $\Psibar$.
By Lang's theorem in $G_M$, we can find $h\in G_M$
such that $h^Fh^{-1}=\wdot$.
On the other hand, $g$ satisfies $g^Fg^{-1}=t\wdot$ for some $t\in T_0$.
Now the map $\wdot F$ is a nonstandard Frobenius endomorphism
since $\wdot^F=\wdot$ and
so $(\wdot F)^m=F^m$, where $m$ is the order of $\wdot$.
Furthermore ${T_0}^{\wdot F}=T_0$.
So, by Lang's theorem in $T_0$, there is a $u\in T_0$ such that
$t=u^{\wdot F}u^{-1}$.
Set $\tilde{g}=u^{-\wdot}g$, so that
$$
  \tilde{g}^F\tilde{g}^{-1}=u^{-\wdot F}g^Fg^{-1}u^{\wdot}
  = u^{-\wdot F}tu\wdot=\wdot
$$
and ${H_0}^{\tilde{g}}={H_0}^{\wdot^{-1}u^{-1}\wdot g}={H_0}^g=H$.
Hence $h^Fh^{-1}=\tilde{g}^F\tilde{g}^{-1}$, that is $\tilde{g}h^{-1}$ is
defined over $k$ and so
$H^{h^{-1}}={H_0}^{\tilde{g}h^{-1}}$ is split.
So
$$
  [M,H^{h^{-1}}]=[M^{h},H]^{h^{-1}}=
  [M,H]^{h^{-1}}\le M^{h^{-1}}=M
$$
and $M$ is split fundamental by Theorem~\ref{T-fund}(\ref{T-fund-splitext}).
\end{proof}

An immediate application is
Algorithm~\ref{A-comp} for computing the direct sum
decomposition of a the Lie algebra $L$.
Although more that one $w$-orbit of $\Phi$ can have the same generalised root,
this clearly is not possible for orbits in different components of $\Phi$.
\begin{figure}
\begin{tabbing}
\quad\=\quad\=\qquad\=\quad\=\kill
$\text{\sc Components} := \text{\bf function}(L)\quad$\\
\> \tbf{let} $H:=\text{\sc MaximalToralSubalgebra}(L)$\\
\> \tbf{let} $\sF=\text{\sc GeneralisedRoots}(L,H)$\\
\> \tbf{for} each $f\in\sF$ \tbf{let} $M_f$ be the subalgebra generated by
$L_f$\\
\> \tbf{construct} the graph with vertices $\sF$\\
\> \> \> and an edge $(f,g)$ whenever $M_f\cap M_g$ is not contained in $H$\\
\> \tbf{let} $C=\emptyset$\\
\> \tbf{for} each graph component $c$ \tbf{do}\\
\> \>  \tbf{add} the subalgebra $\langle M_f\mid f\in c\rangle$ to $C$\\
\> \tbf{end for}\\
\> \tbf{return} $C$\\ 
\textbf{end function}
\end{tabbing}
\caption{Direct sum components}
\label{A-comp}
\end{figure}
The components returned are fundamental
subalgebras.

\subsection{Finding a split maximal toral subalgebra}
Suppose now that we have found a nontrivial split fundamental subalgebra $M$.
The following proposition shows that we can use recursion to find a split
maximal toral subalgebra of $L$.
\begin{proposition}\label{P-cent}
  Suppose that the characteristic of $k$ is greater than $3$.
  Let $M$ be a split fundamental subalgebra of $L$.
  Let $K$ be a split maximal toral subalgebra of $M$.
  Then $C_L(K)$ is a split fundamental subalgebra of $L$ with full rank $n$.
  Hence a split maximal toral subalgebra of $C_L(K)$ is also a split
  maximal toral subalgebra of $L$.
\end{proposition}
\begin{proof}
Let $G_M$ be the split fundamental subgroup of $G$ such that $L(G_M)=M$.
Let $T$ be a split maximal torus of $G$ which normalises $G_M$
and let $H=L(T)$.
Then $U=G_M\cap T$ is a split maximal torus of $G_M$.
By Theorem \ref{T-toral}(\ref{T-toral-splitconj}), we can assume without loss
of generality that $K=L(U)$.

Let $C=C_G(U)$.
Clearly $C$ is normalised by $T$ and it is reductive by 
\cite[Corollary~26.2A]{Humphreys75Book}, hence it is split fundamental.

Let $\Psi$ be the subset of $\Phi=\Phi(G,T)$ consisting of roots of $G_M$, 
or equivalently of $M$.
Let $\Psi'$ be the elements of $\Phi$ which are orthogonal to all elements of
$\Psi$.
Then $\Psi'$ is the root system of $C$ and
$$
  L(C)=H\ds\bigds_{\al\in\Psi'}L_\al.
$$
Clearly $L(C)\le C_L(K)$.

Conversely, suppose $x\in C_L(K)$.
Let $\{e_\al,h_i\}$ be a Chevalley basis of $L$ with respect to $H$ and 
write
$$
  x = \sum_{i=1}^n t_i h_i + \sum_{\al\in\Phi} a_\al e_\al.
$$
If $\al\notin\Psi'$ then there exists $\be\in\Psi$ such that
$\langle\al,\be^\star\rangle\ne0$.
By the basic properties of root data $|\langle\al,\be^\star\rangle|\le3$, so
$\langle\al,\be^\star\rangle$ is still nonzero considered as an element of $k$.
Now $h_\be=[e_{-\be},e_\be]$ is in $K$, and so $[x,h_\be]=0$.
But the coefficient of $e_\al$ in $[x,h_\be]$ is
$a_\al\langle\al,\be^\star\rangle$ by (\ref{E-eh}).
Hence $a_\al=0$ for all $\al\notin\Psi'$ and so $x\in L(C)$.

The second conclusion is an immediate consequence of the first.
\end{proof}

We now have a method for finding split maximal toral subalgebras of $L$:
Find a maximal toral subalgebra $H$, and compute its generalised roots.
For each generalised root $f$, construct the subalgebra
$M_f$ generated by $L_f$.
Now, assuming that we can find $A\subseteq\sF$ for which
$M_A=\langle\sum_{f\in A}M_f\rangle$ is known to be split fundamental and strictly contained in $L$,
find a split maximal toral subalgebra $H_A$ of $M_A$.
By Proposition~\ref{P-cent}, a split maximal toral subalgebra of $C_L(H_A)$ is a
split maximal toral subalgebra of $L$. 
Since $M_A$ and $C_L(H_A)$ are split fundamental
subalgebras of $L$, Theorem~\ref{T-fund}(\ref{T-fund-splitext}) ensures that
they are also the Lie algebras 
of $k$-split connected reductive algebraic groups and so this recursion is valid.
Algorithm~\ref{A-mts} gives the precise method we use.
\begin{figure}
\begin{tabbing}
\quad\=\quad\=\quad\=\quad\=\kill
$\text{\sc SplitMaximalToralSubalgebra} := \text{\bf function}(L,Z)\qquad$ [$Z\le Z(L)$]\\
\> \textbf{repeat}\\
\> \> \textbf{let} $H/Z=\text{\sc MaximalToralSubalgebra}(L/Z)$\\
\> \> \textbf{if} $H$ is split \textbf{then return} $H$\\                                          
\> \> \textbf{let} $\sF= \text{\sc GeneralisedRoots}(L/Z,H/Z)$\\                      
\> \> \textbf{for} each $f\in\sF$ \textbf{compute} $M_f=\phi(\langle(L/Z)_f\rangle)$\\                                               
\> \textbf{until} we find $A\subseteq\sF$ such that $M_A=\langle\sum_{f\in
A}M_f\rangle<L$ 
is split fundamental\\
\> \textbf{let} $H_A=\text{\sc SplitMaximalToralSubalgebra}(M_A,\phi(Z(M_A/Z)))$\\
\> \tbf{let} $C_A=\phi(C_{L/Z}(H_A))$ and $Z=\phi(Z(C_A/Z))$\\
\> \tbf{let} $K=Z$\\
\> \tbf{for} $M$ in $\text{\sc Components}(C_A,Z)$ \tbf{do}\\
\> \> \tbf{let} $K=K+\text{\sc SplitMaximalToralSubalgebra}(M,Z)$\\
\> \tbf{end for}\\
\> \tbf{return} $K$\\
\textbf{end function}
\end{tabbing}
\caption{Finding a split maximal toral subalgebra}
\label{A-mts}
\end{figure}
Note that the second argument $Z$ passed to this function is intended to
indicate that we have a basis of $L(k)$ extending a basis of $Z(k)$, and
the pullback map $\phi:L/Z\to L$.
We take $Z=Z(L)$ initially.
In Section~\ref{S-prob}, we give a method for ensuring that $M_A$ is
known to be split fundamental.

\subsection{Finding a Chevalley basis}
We start by giving a recognition theorem for a standard Chevalley basis.
\begin{theorem}
Suppose the finite field $k$ has characteristicgreater than $3$.
Let $G$ be a $k$-split connected reductive linear algebraic group 
defined over $k$  and let $L$ be the Lie algebra of $G$.
Let $H$ be a $k$-split maximal toral subalgebra of $L$
and let $L=H\bigds_\al L_\al$ be the root decomposition of $L$.
Suppose we have a basis of $L$ consisting of $h_i\in H$ for $i=1,\dots,n$
and $e_\al\in L_\al$ for $\al\in\Phi$.
Further suppose this basis satisfies the equation
$$
  [e_{-\al},e_\al]=\sum_{i=1}^n\langle e_i,\al^\star\rangle h_i
$$ 
for every simple root $\al$
and the equations
$$
  [e_\al,e_\be]=N_{\al\be}e_{\al+\be} \quad\text{and}\quad
  [e_{-\al},e_{-\be}]=N_{-\al,-\be}e_{-\al-\be}
$$
for every extraspecial pair $(\al,\be)$.
Then this is a standard Chevalley basis.
\end{theorem}
\begin{proof}
We need to prove that this basis satisfies the defining equations given in
Subsection~\ref{SS-chevtransf}.
Equation~(\ref{E-hh}) follow from the fact that a toral subalgebra is
abelian,
Equation~(\ref{E-hh}) is given,
and the Equation~(\ref{E-ee}) follows from \cite[Theorem~4.2.1]{Carter72Book}.
It remains to prove Equation~(\ref{E-eh}).

For $y\in Y$, define $h_y=\sum_{i=1}^n\langle e_i,y\rangle h_i\in H(k)$.
It suffices to prove that
\begin{align}
    [e_\alpha,h_y] = \langle\alpha,y\rangle\,e_\alpha,\label{E-y}
\end{align}
for some collection of elements $y$ generating $Y$.

Now~(\ref{E-y}) is true for all $y\in\Phi^\star$ by 
\cite[Theorem~4.2.1]{Carter72Book}.
If $\langle\al,y\rangle=0$ for all $\al\in\Phi$, then $h_y$ is central and so 
(\ref{E-y}) is trivially true.
Together, these two kinds of element generate $Y$ and so we are done.
\end{proof}

A consequence of this theorem is Algorithm~\ref{A-chev} for finding a Chevalley
basis of $L$.
\begin{figure}
\begin{tabbing}
\quad\=\quad\=\quad\=\quad\=\kill
$\text{\sc StandardChevalleyBasis}:=\text{\bf function}(L)$\\
\> \textbf{let} $H=\text{\sc SplitMaximalToralSubalgebra}(L)$\\
\> \textbf{compute} the root system $\Phi$ and root spaces $L_\al$ for $\al\in\Phi$\\
\> \tbf{find} simple roots $\al_1,\dots,\al_\ell$ for $\Phi$\\
\> \tbf{for} $i=1,\dots,\ell$ \tbf{do}\\
\> \> \tbf{let} $\al=\al_i$\\
\> \> \tbf{choose} nonzero $e_{\al}\in L_{\al}$ and $f_{\al}\in L_{-\al}$\\
\> \> \tbf{find} $a\in k$ such that $[e_\al,[f_\al,e_\al]]=2ae_\al$\\
\> \> \tbf{let} $e_{-\al}=f_\al/a$, $h_\al=[e_{-\al},e_\al]$\\
\> \tbf{end for}\\
\> \tbf{compute} a basis $\{h_i\}$ for $H(k)$ with $h_\al=\sum_i\langle
e_i,\al^\star\rangle h_i$ for simple roots $\al$\\
\> \textbf{for} $\ga$ a nonsimple root \tbf{do}\\
\> \> \textbf{let} $(\al,\be)$ be the extraspecial pair of $\ga$\\
\> \> \textbf{let} $e_\ga=[e_\al,e_\be]/N_{\al\be}$,
$e_{-\be}=[e_{-\al},e_{-\be}]/N_{-\al,-\be}$\\
\> \textbf{end for}\\
\> \tbf{return} $\{e_\al,h_i\}$\\
\textbf{end function}
\end{tabbing}
\caption{Finding a standard Chevalley basis}
\label{A-chev}
\end{figure}
Note that for an extraspecial pair $(\al,\be)$, we have $0<N_{\al\be}\le 3$,
so division by $N_{\al\be}$ is not a problem.
The basis $\{h_i\}$ can be computed by elementary linear algebra.
Note that in the second for-loop, the roots are taken in the linear order $<$ of
Subsection~\ref{SS-chevtransf}, thus ensuring that
$e_\al$ and $e_\be$ are already known when we compute $e_\ga$.

\section{Time analysis}\label{S-prob}
Let $L$ be the Lie algebra of the $k$-split connected reductive linear algebraic
group $G$.
We now find bounds on the probability of finding
a maximal toral subalgebra $H\le L$ and a set $A$ of generalised roots such that
$M_A$ is known to be split fundamental.
To simplify our analysis, we just bound the probability that
Algorithm~\ref{A-toral} finds a maximal toral subalgebra in a single step, or
equivalently that the random element chosen is regular semisimple.
Subsection~\ref{SS-rss} gives bounds on the frequencies of 
regular semisimple elements corresponding to Weyl group elements.
In Section~\ref{SS-heur}, we bound the proportion of suitable Weyl group elements.
We give the proof of Theorem~\ref{T-time} in Section~\ref{SS-time}.

Throughout this section, $n$ is the reductive rank of $G$, $\ell$ is the
semisimple rank of $G$, $d$ is the dimension of $L$, and $d_1,\dots,d_\ell$ are the
invariant degrees of $G$ as defined in \cite[Section~9.3]{Carter72Book}.

\subsection{Regular semisimple elements}\label{SS-rss}
An element of $L$ is \emph{regular semisimple} if its centraliser is a maximal
toral subalgebra.
For any subvariety $S$ of $L$, let $S_\rss$ be the variety of regular semisimple
elements in $S$. 
Recall from Subsection~\ref{SS-toral} that the maximal toral subalgebras of $L$ are classified up to
$G(k)$-conjugacy by the conjugacy classes of $W$.
Fix $w$ in $W$ and
let $L_{\rss,w}$ be the set of elements $x\in L$ which are regular semisimple
and such that there exists $g\in G$ with $C_L(x)={H_0}^g$ and 
$g^Fg^{-1}\in T_0\wdot$. 
Although we give direct proofs, many results in this section also follow
from Gus Lehrer's analysis of hyperplane complements \cite{Lehrer92,Lehrer98}.

The following result bounds our chances of finding a regular semisimple element
in $L(k)$ whose centraliser corresponds to the $W$-class of a given $w$.
\begin{proposition}\label{P-prob} 
Let $L$ be the Lie algebra of a $k$-split connected reductive group $G$
with root datum $(X,\Phi,Y,\Phi^\star)$.
Let $w$ be an element of the Weyl group $W$.
Define 
$$
  Q_w(X)= \frac{\prod_{i=1}^\ell(1-X^{d_i})}{\det_Y(1-wX)}.
$$
Then
$$
  \left(1-\sum_{i=1}^\ell\frac{c_i}{q^i}\right) 
  Q_w(1/q) \frac{|w^W|}{|W|}
  \le \frac{|L_{{\rm rss},w}(k)|}{|L(k)|}
  \le Q_w(1/q) \frac{|w^W|}{|W|}.
$$
where
$c_i=c_i(w)$ is the number of $w$-orbits in $\Phi$ consisting of roots $\al$
with the property that $i$ is the largest integer for which 
$\al,\al w,\dots,\al w^{i-1}$ are $\kbar$-linearly independent.
\end{proposition}
\begin{proof}
Fix some $g\in G$ such that $g^Fg^{-1}=\wdot$ and define $H_w={H_0}^g$.
Let $T_w={T_0}^g$ so that $L(T_w)=H_w$.
Then
$$
  L_{\rss,w}(k) =
  \{ x\in L_\rss(k)\mid x\in H_w(k)^h \mbox{ for some } h\in G(k)\},
$$
which is in one-to-one correspondence with
$$
  \{ (x,H)\in L_\rss(k)\times H_w(k)^{G(k)}\mid x\in H\}.
$$
Since $N_{G(k)}(H_w(k))/ T_w(k) \cong C_W(w)$,
we have
$|H_w(k)^{G(k)}| = \frac{|G(k)|}{|T_w(k)| |C_W(w)|}$.
Hence 
$$
  \frac{|L_{\rss,w}(k)|}{|L(k)|} =
  |(H_w)_\rss(k)| \frac{|G(k)|}{|L(k)||T_w(k)|}
  \frac{|w^W|}{|W|} .
$$

Given a root $\al\in\Phi$, define 
$$
  H_{\al} = \{ h\in H_w \mid \al^g(h)=0 \}.
$$
Then $H_{\al}$ is a hyperplane in $H_w$ and 
$(H_w)_\rss=H_w-\bigcup_{\al\in\Phi}H_{\al}$.
Now  ${H_\al}^F=H_{\al w}$,
so $H_\al(k) = \left(\bigcap_j H_{\al w^j}\right)(k)$.
This space has codimension $i$, the largest integer such that $\al,\al
w,\dots,\al w^{i-1}$ are linearly independent.
So, for each $i=1,\dots,\ell$, we are removing $c_i$ subspaces of codimension
$i$ from a $k$-space of dimension $n$.
Hence
$$
  q^n\left(1-\sum_{i=1}^\ell\frac{c_i}{q^i}\right) \le 
  \left|(H_w)_\rss(k)\right| \le q^n.
$$

Using Theorem~9.4.10 of \cite{Carter72Book} and the fact that our group is untwisted,
we get
$$
  |G(k)| = q^{d}\prod_{i=1}^\ell\left(1-\frac1{q^{d_i}}\right).
$$
Using Proposition~3.3.5 of \cite{Carter93Book} and the fact that $F$ is the
standard Frobenius, we find that $T_w(k)$ has order $\det_Y(qI-w)$.
Hence
$$
  \frac{|G(k)|}{|L(k)||T_w(k)|} =  \frac{q^{d}\prod_i(1-1/{q^{d_i}})}{q^{d}\det_Y(qI-w)}
  = \frac{Q_w(1/q)}{q^{n}}.
$$
\end{proof}

The following useful lemma can be proved by elementary calculus.
\begin{lemma}\label{L-bnd}
Let $a_1,\dots,a_m$ be a sequence of nonnegative integers and
suppose that no integer appears more than $a$ times in this sequence.
Then
$$
  \prod_i \left(1-\frac1{q^{a_i}}\right) \ge \left(1-\frac1q\right)^{2a}.
$$
\end{lemma}

\subsection{Reflection derangements}\label{SS-heur}
Recall from Subsection~\ref{SS-rtdec} that there is a relationship between
between the generalised roots $f$ with respect to a toral subalgebra
and the orbits of the corresponding Weyl group element $w$ on $\Phi$.
This relationship need not be a one-to-one correspondence.  
As we saw in Lemma~\ref{L-orbit}(\ref{L-orbit-fix}) and~(\ref{L-orbit-2}),
this relationship is almost a one-to-one correspondence when the
degree of $f$ is one, or the degree is two and $f=f_-$.
This happens when there is a root $\al$ such that $\al w=\pm\al$.
In other words, when a reflection $s_\al$ is fixed under conjugation by $w$.

In this section, we count the number of Weyl group elements of this kind. 
Given a permutation representation of a group, an element of
the group is called a \emph{derangement} with respect to the representation
if it fixes no points at all. 
The proportion of derangements of the symmetric
group $\Sym_m$ acting on $m$ letters is known to approach $1/e$ as $m\to\infty$.
We give similar results for a Weyl group
acting on its reflections by conjugation.
We refer to these elements as \emph{reflection derangements}.
We are grateful to Anthony Henderson for helping us with the proof of this
proposition.
\begin{proposition}\label{P-e}
If $W$ is an irreducible Coxeter group of classical type $\CA_\ell$, $\CB_\ell/\CC_\ell$, or
$\CD_\ell$, then the proportion of its reflection
derangements approaches $2e^{-3/2}$,
$e^{-5/4}$, $2e^{-5/4}+(4e)^{-1}$, respectively, as $\ell\to\infty$.
For exceptional types, the proportions are as  listed below:\smallskip\\
\centerline{
\begin{tabular}{lllll}\hline
$\CG_2$ & $\CF_4$ & $\CE_6$     & $\CE_7$     & $\CE_8$ \\
$1/3$   & $1/4$   & $1409/2592$ & $1646/2835$ & $3385549/6220800$\\\hline
\end{tabular}}
\end{proposition}

\begin{proof}
Denote by $f$ the number of reflection derangements of $W$.
We wish to determine $f/|W|$.  
\medskip 

\noindent{\bf Type $\CA_\ell$:}
The Weyl group $W(\CA_\ell)$ can be identified with the symmetric group
$\Sym_{\ell+1}$ on $\ell+1$ letters.  
Write $m=\ell+1$ and write $d_m$ for the proportion of
permutations in $\Sym_{m}$ without fixed points in $\{1,\ldots,m\}$.
Denote by $R_m$
the set of all permutations in $\Sym_{m}$ with at most one
fixed point in $\{1,\ldots,m\}$.

An element of $\Sym_m$ does not fix a reflection if, and only if,
it belongs to $R_m$ and does not contain
a transposition $(i,j)$ in its cycle decomposition.
So
$$f = \left|R_m- \bigcup_{1\le i< j\le m}R_m^{i,j}\right|$$
where
$$R_m^{ij} = \{w\in R_m\mid w \mbox{ contains }(i,j) \}.
$$
We compute $f$ by inclusion/exclusion.  As
$R_m^{ij}$ and $R_m^{ij'}$ intersect trivially 
for $j\ne j'$ we can find $f$ as an
alternating sum
over $h$-tuples of commuting transpositions:
$$\sum_{h=0}^{\lfloor m/2\rfloor} (-1)^h\binom{m}{2h}
\frac{(2h)!}{2^h h!} |R_{m-2h}|.$$
Since,
clearly, $|R_m| = d_m + m d_{m-1}/m$,
$$f =  m! \sum_{h=0}^{\lfloor m/2\rfloor} \left(-\frac{1}{2}\right)^h
\frac{1}{h!} \left(d_{m-2h} + d_{m-2h-1}\right)
.$$
As $\lim_{m\to\infty}d_m = 1/e$,
the required proportion
tends to
$$\lim_{m\to\infty} \frac{f}{m!} =
 \sum_{h=0}^{\infty} \left(-\frac{1}{2}\right)^h
\frac{1}{h!} \frac{2}{e} = e^{-\frac{1}{2}}\frac{2}{e} = 2 e^{-\frac{3}{2}}
.$$\medskip 

\noindent{\bf Types $\CB_\ell$ and $\CC_\ell$:}
The Weyl group $W = W(\CB_\ell)=W(\CC_\ell)$ can be identified with the group of
all permutations $w$ of $\{\pm1,\dots,\pm\ell\}$ such that
$(-i)w=-(iw)$.
Define the homomorphism $\phi: W \to\Sym_\ell$ by $iw^\phi=|iw|$.
Then $w\in W$ fixes no reflections if, and only if,
$w^\phi$ is a derangement of $\Sym_\ell$
and, for every transposition $(i,j)$ contained in
the cycle decomposition of $w^\phi$, either
$(i,j,-i,-j)$ or $(j,i,-j,-i)$ is contained in the cycle decomposition of $w$.

Writing $S_\ell$ for elements of $W$ such that $w^\phi$ is a
derangement and 
$$
  S_\ell^{ij} = \{w\in S_\ell \mid w \mbox{ contains } (i,j)(-i,-j) 
  \mbox{ or } (i,-j)(-i,j)\},
$$ 
we find that 
$$
f=\left|S_\ell- \bigcup_{1\le i< j\le \ell}S_\ell^{ij}\right|.
$$
Again,
we can count $f$ by
taking alternating sums
over $h$-tuples of commuting transpositions in $W^\phi$.
As each transposition in the decomposition of an element
of $w^\phi$ corresponds to two 4-cycles as indicated above, we find an
extra factor $2^h$ compared to the $A_\ell$ case:
$$
  \sum_{h\ge0, 2h\le \ell} (-1)^h\binom{\ell}{2h}
  \frac{(2h)!}{2^h h!} 2^h |S_{\ell-2h}|.
$$
As $|S_\ell| = 2^\ell\, \ell!\, d_\ell$,
$$
  f = \sum_{h\ge0, 2h\le \ell} \left(-1\right)^h
  \frac{\ell!}{h!}2^{\ell - 2h}d_{\ell-2h}.
$$
As $\lim_{m\to\infty}d_m = 1/e$
and $|W(\CB_\ell)| = 2^\ell \ell!$,
the required proportion
tends to
$$\lim_{m\to\infty} \frac{f}{2^\ell\ell!} =
 \sum_{h=0}^{\infty} \left(-\frac{1}{4}\right)^h
\frac{1}{h!} \frac{1}{e} = e^{-\frac{1}{4}}e^{-1} =  e^{-\frac{5}{4}}
.$$\medskip 

\noindent{\bf Type $\CD_\ell$:}
The Weyl group $W(\CD_\ell)$ is the subgroup of $W(\CB_\ell)$ consisting of all
elements $w$ such that $\prod_{i=1}^\ell iw$ is positive.
In cycle notation, this means that $w$ has an even number
of negative cycles (that is, cycles in which both positive and negative
numbers occur).

Define $\phi: W \to\Sym_\ell$ as the restriction of the map for type $\CB_\ell$.
Then $w\in W$ does not commute with any reflection if, and only if,
\begin{enumerate}[(i)]
\item
$w^\phi$ fixes at most one element of $\{1,\ldots,\ell\}$
and, for every transposition $(i,j)$ contained in
the cycle decomposition of $\phi(w)$, the cycle occurring in $w$ is
$(i,j,-i,-j)$ or $(j,i,-j,-i)$; or
\item $w^\phi$ has exactly two fixed points, say $i$ and $j$,
and the cycle decomposition of $w$ contains $(i,-i)(j)(-j)$ or $(i)(-i)(j,-j)$.
\end{enumerate}
The number of elements of the type (ii)
is clearly $\binom{\ell}{2} d_{\ell-2} 2^{\ell -2}(\ell-2)!$,
contributing
$$
  \lim_{\ell\to\infty}
  \frac{\binom{\ell}{2}  2^{\ell -2} d_{\ell-2}(\ell-2)!}{ |W(\CD_\ell)|}
  = \lim_{\ell\to\infty} 2^{-2} {d_{\ell-2}} = \frac{1}{4e}
$$
to the required asymptotic proportion.

Writing $T_\ell$
for elements of $W$ such that $w^\phi$ fixes at most two
elements and
$$
  T_\ell^{i,j} = \{w\in T_\ell\mid w \mbox{ contains }
  (i,j)(-i,-j) \mbox{ or } (i,-j)(j,-i)\},
$$
we find
that the set of elements of type (i) is
$$
  T_\ell - \bigcup_{1\le i< j\le \ell}T_\ell^{i,j}.
$$
Again,
we take alternating sums
over $h$-tuples of commuting transpositions in $\phi(W)$.
As each transposition in the decomposition of an element
of $\phi(w)$ corresponds to two 4-cycles as indicated above, we find the same
factor $2^h$ as for the $B_\ell$ case:
$$
  \sum_{h=0}^{\lfloor\ell/2\rfloor} (-1)^h\binom{\ell}{2h}
  \frac{(2h)!}{2^h h!} 2^h |T_{\ell-2h}|.
$$
As $|T_\ell| = 2^{\ell-1} \,\ell!\,( d_\ell +  d_{\ell - 1})$,
the result is
$$ 
  \sum_{h=0}^{\lfloor\ell/2\rfloor} \left(-\frac{1}{4}\right)^h
  \left(d_{\ell-2h}+d_{\ell-2h-1}\right),$$
which contributes
$$\lim_{m\to\infty} \frac{f}{2^\ell\ell!} =
 \sum_{h=0}^{\infty} \left(-\frac{1}{4}\right)^h
\frac{1}{h!} \frac{2}{e} = 2e^{-\frac{1}{4}}e^{-1} = 2 e^{-\frac{5}{4}}
$$
to the required proportion.
Hence, the asymptotic proportion is
$(4e)^{-1} +  2 e^{-{5}/{4}}$.\medskip 

\noindent{\bf The exceptional types:} These were computed by machine.
\end{proof}

\begin{corollary}\label{C-e}
The proportion of reflection derangements in a Weyl group is less than
$\frac{2}{3}$.
\end{corollary}
\begin{proof}
Recall that if $a_n>0$ converges monotonically to zero, then 
$\sum_{i=0}^\infty(-1)^ia_i$ is
called an \emph{alternating series}.
The maximum value of the partial sums 
$s_n=\sum_{i=0}^n(-1)^ia_i$
of such a series is one of the first two
partial sums.
Since the series in the previous proposition are sums of alternating sequences,
it is always possible to find a constant $M$ such that the maximum value of the
partial sums is one of $s_1,\dots,s_M$.
It is now easy to show on a case-by-case basis that
the proportion of reflection derangements in an irreducible Weyl group is at
most $\frac{2}{3}$.

If $W$ is a direct product decomposition into $s$
irreducible Weyl groups, then an element of $W$ is a reflection derangement
if and only if each component of $w$ is
a reflection arrangement, and so their proportion is at
most $\left(\frac{2}{3}\right)^s\le\frac23$.
\end{proof}
Together with Proposition~\ref{P-prob}, this shows that the chance of finding a
regular semisimple element of $L$ corresponding to a reflection nonderangement
in the Weyl group is at least one third, provided $q$ is large enough.
To complete the analysis, we need a more precise bound on the probability of
finding certain regular semisimple elements.

\subsection{Time analysis}\label{SS-time}
We start by looking at the Coxeter class in the Weyl group.
The Coxeter element is actually a reflection derangement, but this proof
is the model for our next result.
\begin{proposition}\label{P-cox}
Suppose that $W$ is an irreducible Weyl group.
If $w_c$ is a Coxeter element of $W$, then
$$
  \frac{|L_{{\rm rss},w_c}(k)|}{|L(k)|}\ge
  \left(1-\frac\ell{q^{\ell/2}}\right)
  \left(1-\frac{1}{q}\right)^4\frac1{h}
$$
where $h$ is the order of $w_c$.
\end{proposition}
\begin{proof}
Suppose $\al$ is a root and $\al {w_c}^m$ is a linear combination of
$\al,\al w_c,\dots,\al {w_c}^{m-1}$.
We prove that $m\ge\ell/2$ on a case-by-case basis:\medskip 

\noindent{\bf Type $\CA_\ell$:} Identify $W$ with $\Sym_{\ell+1}$ and consider $\Phi$
  to consist of roots $e_i-e_j$ with $i\ne j$.
  We can take $w_c=(1,2,\dots,\ell+1)$ and $\al=e_i-e_j$.
  So $\al {w_c}^m=e_{i+m}-e_{j+m}$ with the subscripts taken modulo $\ell+1$.
  Hence $\al {w_c}^m$ is a linear combination of $\al,\dots,\al {w_c}^m$ if{f}
  $i+m$ and $j+m$ are both in $[i,i+m-1]\cup[j,j+m-1]$ modulo $\ell+1$.  
  By the pigeon hole principle, this can only happen if $m\ge(\ell+1)/2$.\medskip 

\noindent{\bf Type $\CC_\ell$:}
  Identify $W$ with the set of permutations $w$ of $\{\pm1,\dots,\pm\ell\}$
  such that $(-i)w=-(iw)$ for $i=1,\dots,\ell$.
  Consider $\Phi=\Phi(C_\ell)$ to consist of roots $\ep e_i-\del e_j$ with
  $\ep,\del\in\{\pm1\}$, $i,j=1,\dots,\ell$ and $\ep i\ne\del j$.
  We can take $$w_c=(1,2,\dots,\ell,-1,-2,\dots,-\ell)$$ and 
  $\al=\ep e_i-\del e_j$.  
  The same argument used in type $\CA_\ell$ now shows that $m\ge\ell/2$. \medskip 

\noindent{\bf Type $\CB_\ell$:}
  The permutation action of $W(\CB_\ell)$ on its roots is isomorphic to the
  action of $W(\CC_\ell)$ on its roots, so the same argument works.\medskip 

\noindent{\bf Types $\CD_\ell$:}
  Identify $W$ with the elements of $W(\CC_\ell)$ such that 
  $\prod_{i=1}^\ell (iw)>0$ and consider $\Phi$ to consist of the roots 
  $\ep e_i-\del e_j$  with $\ep,\del\in\{\pm1\}$, $i,j=1,\dots,\ell$ and $i\ne j$.
  We can take $w_c=(1,2,\dots,\ell-1,-1,-2,\dots,-\ell+1)(\ell,-\ell)$
  and $\al=\ep e_i-\del e_j$.
  Once again $m\ge\ell/2$ if $i,j\ne\ell$.
  If $i=\ell$, $j\ne\ell$, then $\al {w_c}^m=(-1)^m\ep e_{\ell}-\del e_{j+m}$
  with the second subscript taken modulo $\ell-1$, and so $m\ge\ell-1$.\medskip 

\noindent{\bf Exceptional types:} These are easily checked by computer. \medskip 

It is well known that every orbit of $w_c$ on $\Phi$ has size $h$, so
$\sum_ic_i(w_c)=2N/h=\ell$.
We have shown that $c_i(w_c)=0$ for $i<\ell/2$, so
$$
  1-\sum_{i=1}^\ell\frac{c_i(w_c)}{q^i} \ge 1-\frac\ell{q^{\ell/2}}.
$$

The functions $Q_{w_c}(X)$ are straightforward to
compute and are given in Table~\ref{T-pols}.
\begin{table}
\begin{tabular}{ll} \hline
$\CA_\ell$ & $\prod_{i=1}^\ell(1-X^i)$\\
$\CB_\ell,\CC_\ell$ & $(1-X^\ell)\prod_{i=1}^{\ell-1}(1-X^{2i})$\\
$\CD_\ell$ & $\prod_{i\in\{1,\ell-1,\ell\}}(1-X^i)
            \prod_{i=2}^{\ell-2}(1-X^{2i})$\\
$\CG_2$    & $(1-X^2)(1-X^3)(1+X)$\\
$\CF_4$    & $(1-X^6)\prod_{i\in\{4,6,8\}}(1-X^i)$\\
$\CE_6$    & $(1-X^6)\prod_{i\in\{1,4,5,6,8\}}(1-X^i)(1+X^3+X^6)$\\
$\CE_7$    & $(1-X^6)\prod_{i\in\{1,6,8,10,12,14\}}(1-X^i)(1+X^3+X^6)$\\
$\CE_8$    & $\prod_{i\in\{1,8,10,12,14,18,20,24\}}(1-X^i)
            (1+X^3)(1+X^5+X^{10})$\smallskip\\\hline\smallskip
\end{tabular}
\caption{The functions $Q_w(X)$ for a Coxeter element $w$}\label{T-pols}
\end{table}
The terms in which every coefficient is positive can be ignored, since they are
bounded below by 1 when we set $X=1/q$.
Since no term $1-X^a$ appears more than twice in these polynomials
and $q\ge 3$,
it follows by Lemma~\ref{L-bnd} that $Q_{w_c}(1/q)\ge (1-1/q)^4$.

The required inequality now follows from the first inequality of
Proposition~\ref{P-prob} and the fact that the centraliser of $w_c$ has order
$h$.
\end{proof}

We now consider reflection nonderangements that are, in some sense, close to
being Coxeter elements.
\begin{proposition}\label{P-subcox}
Suppose that $W$ is an irreducible Weyl group of rank greater than one.
If $W$ is classical with rank at least $7$ then
there is a reflection nonderangement $w$ such that
$$
  \frac{|L_{{\rm rss},w}(k)|}{|L(k)|} \ge
  \left(1-\frac3q-\frac4{q^2}-\frac{\ell+5}{q^{(\ell-2)/2}}\right)
  \left(1-\frac{1}{q}\right)^6
  \frac1{4\ell}.
$$
For other Cartan types there is a reflection nonderangement $w$ such that
$$
  \frac{|L_{{\rm rss},w}(k)|}{|L(k)|} \ge
  \left(1-\sum_{i=1}\frac{c_i}{q^i}\right)
  \left(1-\frac{1}{q}\right)^6
  \frac1{c}.
$$
with the constants $c$ and $c_i$ listed in Table~\ref{T-cis}.
\end{proposition}
\begin{table}
\begin{tabular}{lllrrrrrrrr}\hline
Type    & $c$ & $c_1$ & $c_2$ & $c_3$ & $c_4$ & $c_5$ & $c_6$ & $c_7$ & $c_8$ \\\hline
$\CA_1$ & 2 & 1 \\
$\CA_2$ & 2 & 1 & 2 \\
$\CB_2$ & 8 & 4 & 0 \\
$\CG_2$ & 4 & 3 & 4 \\
$\CA_3$ & 8 & 2 & 4 & 0 \\
$\CB_3$ & 8 & 1 & 2 & 2 \\
$\CA_4$ & 6 & 1 & 2 & 0 & 2 \\
$\CB_4$ & 12 & 1 & 1 & 4 & 0 \\
$\CD_4$ & 16 & 5 & 8 & 0 & 0 \\
$\CF_4$ & 36 & 3 & 1 & 6 & 0 \\
$\CA_5$ & 8 & 1 & 1 & 2 & 4 & 0 \\
$\CB_5$ & 16 & 1 & 0 & 0 & 4 & 2 \\
$\CD_5$ & 16 & 3 & 5 & 2 & 4 & 0 \\
$\CA_6$ & 10 & 1 & 0 & 0 & 4 & 0 & 2 \\
$\CB_6$ & 20 & 1 & 0 & 0 & 2 & 5 & 0 \\
$\CD_6$ & 24 & 3 & 5 & 3 & 4 & 0 & 0 \\
$\CE_6$ & 36 & 3 & 0 & 3 & 2 & 6 & 0 \\
$\CE_7$ & 60 & 3 & 0 & 0 & 2 & 10 & 0 & 0 \\
$\CE_8$ & 108 & 3 & 0 & 0 & 0 & 0 & 4 & 9 & 0 \\\hline\smallskip
\end{tabular}
\caption{The constants $c$ and $c_i$ for small rank and exceptionals}\label{T-cis}
\end{table}
\begin{proof}
Fix a root $\be$.
Assume $\be$ is short (resp.\ long) for Cartan type $\CB_\ell$ (resp.\
$\CC_\ell$). 
Let
$
  \Phi_\be = \{\ga\in\Phi \mid \langle\ga,\be^\star\rangle=0\}.
$
Then $\Phi_\be$ is a subsystem of $W$ and, except in type $\CD_4$, it
has at most two irreducible components.
Let $\Phi_\be'$ be the irreducible summand of $\Phi_\be$ of maximal rank.
Let $s_\be$ be the reflection in $\be$ and let $w_\be$ be the Coxeter element 
of $W(\Phi_\be')$.
We take $w=s_\be w_\be$, except for type $\CA_1$ where we use $w=1$,
type $\CG_2$ where we use $w=s_\be$, and
type $\CD_4$ where we use $s_1s_2s_1s_3s_2s_1s_4s_2s_1s_3s_2$.
(Here $s_i$ is the $i$th simple reflection, with the numbering given in
\cite{Bourbaki75Book}.)
These elements are all reflection nonderangements.

First we prove that
$$
  \sum_{i=1}\frac{c_i}{q^i} \le \frac3q+\frac4{q^2}+\frac{\ell+5}{q^{(\ell-2)/2}}
$$ 
for the classical types of rank at least 7.
\medskip

\noindent{\bf Type $\CA_\ell$:} Assume $\be=e_1-e_2$.
  Then $\Phi_\be$ has type $\CA_{\ell-2}$, and so orbits within $\Phi_\be$ contribute
  at most $\frac{\ell-2}{q^{(\ell-2)/2}}$ to the sum, as in the previous proof.
  If $\al\notin\Phi_\be$, then $\al=\pm(e_i-e_j)$ where $i=1$ or $2$ and $j>2$.
  These roots form one orbit of size $2$ and either two orbits of size $\ell-1$ or one orbit of size $2(\ell-1)$.
  So these orbits contribute at most $1/q+2/q^{\ell}$.\medskip 

\noindent{\bf Type $\CB_\ell$ with $\be$ short:}
  Assume $\be=e_1-e_2$.
  Then $\Phi_\be$ has type $\CB_{\ell-1}$, and so the orbits within $\Phi_\be$ contribute
  at most $\frac{\ell-1}{q^{(\ell-1)/2}}$.
  If $\al\notin\Phi_\be$, then $\al=\ep e_i-\del e_j$ where $i=1$ or $2$ and $j>2$.
  These roots form four orbits of size two with $m=1$ and four orbits with $m=\ell-2$.\medskip 

\noindent{\bf Type $\CC_\ell$ with $\be$ long:} This is similar to type
$\CB_\ell$, with the short roots and long roots exchanged.\medskip 

\noindent{\bf Type $\CD_\ell$:}
  Assume $\be=e_1-e_2$.
  Then $\Phi_\be$ has type $\CD_{\ell-2}\CA_1$ and
  $\Phi_\be'$ is the subsystem of type $\CD_{\ell-2}$.
  So the orbits within $\Phi_\be'$ contribute
  at most $\frac{\ell-2}{q^{(\ell-1)/2}}$ to the sum.  
  If $\al\notin\Phi_\be'$, then $\al=\ep e_i-\del e_j$ where $i=1$ or $2$ and
  $j>2$. 
  These roots form at most four orbits with $m=\ell-2$.\medskip 

The values of the constants in Table~\ref{T-cis} are easily computed in Magma.
The constant $c$ is just $|C_W(w)|$.
The functions $Q_w(X)$ are given in Table~\ref{T-pols2}.
\begin{table}
\begin{tabular}{ll}\hline
$\CA_1$  & $1-X$\\
$\CA_2$  & $1-X^3$\\
$\CA_3$  & $(1-X)(1-X^3)(1+X^2)$\\
$\CA_\ell$ $(\ell>3)$ & $\prod_{i\in\{1,\dots,\ell+1\}\setminus\{2,\ell-1\}}(1-X^i)$\\
$\CB_2,\CC_2$ & $(1-X)^2(1+X^2)$ \\
$\CB_3,\CC_3$ & $(1-X)(1-X^2)(1-X^6)$ \\
$\CB_4,\CC_4$ & $(1-X)(1-X^3)(1-X^4)(1-X^8)$\\
$\CB_\ell,\CC_\ell$ $(\ell>4$)
   & $(1-X)(1-X^{\ell-2})^{(3-(-1)^\ell)/2}
      \prod_{i\in\{2,\dots,\ell\}\setminus\{\ell-1\}}(1-X^{2i})$\\
$\CD_4$  & $(1-X)^2(1+X^6)(1+X^2)^2$\\
$\CD_5$  & $(1-X)^3(1-X^5)(1-X^6)(1+X^2)(1+X^4)$\\
$\CD_6$  & $(1-X)^3(1-X^3)(1-X^6)(1-X^{10})(1+X^2)^2(1+X^4)$\\
$\CD_\ell$  $(\ell>6$) & $(1-X)^3\prod_{i\in\{\ell-3,\ell\}}(1-X^i)
            \prod_{i\in\{4,\dots,\ell-1\}\setminus\{\ell-3\}}(1-X^{2i})
            \times$\\
           & $\quad(1+X^2)(1+X^2+X^4)$\\
$\CG_2$    & $1-X^6$\\
$\CF_4$    & $\prod_{i\in\{1,3,8,12\}}(1-X^i)$\\
$\CE_6$    & $(1-X)^2\prod_{i\in\{5,8,9,12\}}(1-X^i)$\\
$\CE_7$    & $(1-X)^2\prod_{i\in\{5,6,12,14,18\}}(1-X^i)(1+X^2)(1+X^4)$\\
$\CE_8$    & $(1-X)^2\prod_{i\in\{6,12,14,20,24,30\}}(1-X^i)\times$\\
           & $\quad(1+X^2)(1+X^4)(1+X^3+X^6)$
\smallskip\\\hline\smallskip
\end{tabular}
\caption{The functions $Q_w(X)$}\label{T-pols2}
\end{table}
Applying Lemma~\ref{L-bnd}, we get $Q_{w}(1/q)\ge (1-1/q)^6$.

For groups not covered in Table~\ref{T-cis},
let $h_\be$ be the Coxeter number of $\Phi_\be'$.
Then the centraliser of $w_\be$ in $W(\Phi_\be')$ has order
$h_\be$, and the centraliser of $w$  in $W$ has order
$2h_\be \le  4\ell$.
The required result now follows from the first inequality of
Proposition~\ref{P-prob}.
\end{proof}

Finally we are in a position to give an analysis of our algorithm.
We refer to Algorithm~\ref{A-mtsA1}, a version of
Algorithm~\ref{A-mts} which searches for maximal toral subalgebras corresponding to
reflection nonderangements.
\begin{figure}
\begin{tabbing}
\quad\=\quad\=\quad\=\quad\=\quad\quad\quad\quad\quad\quad\quad\quad\quad\quad\quad\quad\quad\quad\quad\quad\quad\quad\quad\quad\quad\quad\quad\=\kill
$\text{\sc SplitMaximalToralSubalgebra} := \text{\bf function}(L,Z)\qquad$\\
\> \textbf{repeat}\\
\> \> \textbf{let} $H/Z=\text{\sc MaximalToralSubalgebra}(L/Z)$        \>\>\>$O(\ell^6\log(q))$\\
\> \> \textbf{if} $H$ is split \textbf{then return} $H$                \>\>\>$O(\ell^7\log(q))$\\
\> \> \textbf{let} $\sF= \text{\sc GeneralisedRoots}(L/Z,H/Z)$         \>\>\>$O(\ell^7\log(\ell)\log^2(q))$\\
\> \> \tbf{if} there exists $f\in\sF$ with $\deg(f)=1$ \tbf{then}\\
\> \> \> \tbf{let} $M/Z := \langle\phi((L/Z)_f + (L/Z)_{f_-})\rangle$    \>\>$O(\ell^6\log(q))$\\
\> \> \tbf{elif} there exists $f\in\sF$ with $\deg(f)=2$ and $f=f_-$  \tbf{then}\\
\> \> \> \tbf{let} $M/Z := \langle\phi((L/Z)_f\rangle$                   \>\>$O(\ell^6\log(q))$\\
\> \> \> \tbf{if} $\dim(M/Z)\ne3$ \tbf{then}\\
\> \> \> \> \tbf{find} $\al$ in $\Phi_f$ over $k_2$                        \>$O(\ell\log^2(q))$\\
\> \> \> \> \tbf{let} $M/Z=\langle\phi((L/Z)_\al+(L/Z)_{-\al})\rangle$     \>$O(\ell^6\log(q))$\\
\> \> \tbf{end if}\\
\> \textbf{until} $M$ is defined                                     \>\>\>\>$O(\ell\log(\ell))$ times\\
\> \textbf{let} $H=\text{\sc SplitMaximalToralSubalgebra}(M,\phi(Z(M/Z)))$\\
\> \tbf{let} $C=\phi(C_{L/Z}(H))$ and $Z=\phi(Z(C/Z))$\\
\> \tbf{let} $K=Z$\\
\> \tbf{for} $M$ in $\text{\sc Components}(C,H)$ \tbf{do}\\
\> \> \tbf{let} $K=K+\text{\sc SplitMaximalToralSubalgebra}(M,Z)$\\
\> \tbf{end for}\\
\> \tbf{return} $K$\\
\textbf{end function}
\end{tabbing}
\caption{Finding a split maximal toral subalgebra}
\label{A-mtsA1}
\end{figure}
As discussed in Subsection~\ref{SS-heur}, finding $f$ with $\deg(f)=1$, 
or $\deg(f)=2$ and $f=f_-$ is equivalent to the corresponding Weyl group element
being a reflection nonderangement.
When $\deg(f)=2$ and $f=f_-$, we have always found in practice that $M/Z$ is of type
$\CA_1$,
and so has dimension $3$.
We do not have a proof of this however, so it is necessary to check and then decompose over the
field extension $k_2$ in the unlikely event that we get a larger subalgebra.

\begin{theorem}\label{T-mtatime}
Suppose that the characteristic of $k$ is greater than $3$.
Let $G$ be a $k$-split connected reductive group and 
let $L$ be the Lie algebra of $G$.
We can find a split maximal toral subalgebra of $L$ in Las Vegas time
$O(n^3\ell^6\log^2(\ell)\log^2(q))$.
\end{theorem}
\begin{proof}
Before calling Algorithm~\ref{A-mtsA1}, we compute the centre of $L$, which
takes time $O((n+\ell^2)^3\log(q))$. 
Using Algorithm~\ref{A-comp}, we can assume $G$ is simple.
As indicated in Algorithm~\ref{A-mtsA1}, the computations within the main loop take time
$O(\ell^7\log(\ell)\log^2(q))$. 

By Proposition~\ref{P-subcox}, if $G$ is classical with rank at least $7$, 
we obtain a split toral subalgebra $M$ with
probability at least
$$
  \left(1-\frac1q\right)^6
  \left(1-\frac3q-\frac4{q^2}-\frac{\ell+5}{q^{(\ell-2)/2}}\right)
  \frac1{4\ell}.
$$
For $q\ge5$ and $\ell\ge7$, this is at least
$$
  \left(\frac45\right)^6
  \left(1-\frac35-\frac4{25}-\frac{12}{5^{5/2}}\right)\frac1{4\ell}>0.
$$
Similarly for the Cartan types in Table~\ref{T-cis}, except for $\CD_4$,
$$
  \left(1-\sum_{i=1}\frac{c_i}{q^i}\right)
  \left(1-\frac{1}{q}\right)^6
  \frac1{c} \ge   \left(1-\sum_{i=1}\frac{c_i}{5^i}\right)
  \left(1-\frac{1}{5}\right)^6
  \frac1{c}>0.
$$
For type $\CD_4$, the bound is negative for $q=5$, but positive for $q\ge7$.
So it remains to consider the Lie algebra $\CD_4(5)$.
But for any fixed Lie algebra, it is easily
seen that there is a nonzero chance of the algorithm working, since there is a
chance that the toral subalgebra found by Algorithm~\ref{A-toral} is already 
split.
We have now shown that there is a constant $C>0$ such that the probability of
success after one iteration of the main loop is at least $C/\ell$.

Since
$$
  \lim_{\ell\to\infty} \left(1-\frac{C}\ell\right)^{a\ell} = e^{-aC},
$$
we can choose $a$ such that 
$$
  \left(1-\frac{C}{\ell}\right)^{a\ell} \le \frac1{e^4}
$$
for all $\ell$.
Hence the probability of failure after $a\ell\log(\ell)$ repetitions of the 
loop is at most
$$
  \left(1-\frac{C}{\ell}\right)^{a\ell\log(\ell)} \le
  \left(\frac1{e^4}\right)^{\log(\ell)}= \frac1{\ell^4}.
$$

Clearly the depth of recursion is at most $\ell$,
which contributes a further factor of $\ell$ to our timing.
The ranks of
all the subalgebras in all the calls at a particular depth sum to at most
$\ell$, so the total number of recursive calls is at most $\ell^2$.
Hence
the overall probability of success is at least
$$
  \left(1-\frac1{\ell^4}\right)^{\ell^2}\ge 
  \left(1-\frac1{2\ell^2}\right)^{\ell^2} \ge \frac12.
$$ 
Hence Algorithm~\ref{A-mtsA1}
takes Las Vegas time $O(\ell^9\log^2(\ell)\log^2(q))$.
Combining this with the  preprocessing time of $O((n+\ell^2)^3\log(q))$, and using
the fact that $n\ge\ell$ we get the desired result.
\end{proof}

\begin{corollary}\label{C-chevtime}
Suppose that the characteristic of $k$ is greater than $3$.
Let $G$ be a $k$-split connected reductive group and 
let $L$ be the Lie algebra of $G$.
We can find a Chevalley basis of $L$ in Las Vegas time
$O(n^3\ell^6\log^2(\ell)\log^2(q))$.
\end{corollary}
\begin{proof}
The time taken to find a split maximal toral subalgebra clearly dominates
the time for Algorithm~\ref{A-chev}.
\end{proof}
We can easily decompose $G$ into simple subgroups, since we know its root datum.
Hence, combining this corollary with Proposition~\ref{P-ad},
we see that the algorithm for Lang's Theorem takes Las Vegas time
$$
  O(n^3\ell^6\log^2(\ell)\log^2(q)+n^8r^2s^2\log^2(q)),
$$ 
which is easily
simplified to the expression in Theorem~\ref{T-time}.

\bibliographystyle{amsalpha}

\def\cprime{$'$} \def\Dbar{\leavevmode\lower.6gex\hbox to 0pt{\hskip-.23ex
  \accent16\hss}D}
\providecommand{\bysame}{\leavevmode\hbox to3em{\hrulefill}\thinspace}
\providecommand{\MR}{\relax\ifhmode\unskip\space\fi MR }
\providecommand{\MRhref}[2]{%
  \href{http://www.ams.org/mathscinet-getitem?mr=#1}{#2}
}
\providecommand{\href}[2]{#2}

\end{document}